\documentclass[11pt,english,reqno]{amsart}
\usepackage[T1]{fontenc}
\usepackage{mathrsfs}
\usepackage{amsmath,amsthm,amssymb}
\usepackage{latexsym}
\usepackage{times}
\usepackage{enumerate}
\usepackage{mathrsfs}
\usepackage{stmaryrd}
\usepackage{amsopn}
\usepackage{amsmath}
\usepackage{amssymb}
\usepackage{amsfonts}
\usepackage{amsbsy}
\usepackage{amscd,indentfirst,epsfig}
\usepackage{hyperref}
\usepackage{amsfonts,amsmath,latexsym,amssymb,verbatim,amsbsy}
\usepackage{amsthm}
\usepackage{dsfont}
\usepackage{colordvi}
\usepackage{pstricks}
\newtheorem{theo}{Theorem}[section]
\newtheorem{lemme}[theo]{Lemma}

\newtheorem{propo}[theo]{Proposition}
\newtheorem{cor}[theo]{Corollary}
\newtheorem{hyp}[theo]{Assumption}
\newtheorem{nb}[theo]{Remark}

\theoremstyle{definition}

\def \leq {\leqslant}
\def \geq {\geqslant}
\numberwithin{equation}{section}
\def\ind#1{\lower5pt\hbox{$\scriptstyle #1$}}
\def \le {\leqslant}
\def \ge {\geqslant}

\def \d {\, \mathrm{d} }
\def \G {\mathcal{G}}

\def \e {e}

\def \ds {\displaystyle}

\def \a {\alpha}
\def \b {\beta}

\def \g {\overline{\gamma}}
\def \ds {\displaystyle}
\def\Q {\mathcal{Q}}

\def\R{{\mathbb R}}
\def \S {{\mathbb S}^2}
\def \q {q}
\def \n {\mathrm{n}}

\def \u {\mathbf{u}}
\def \v {{v}}
\def \vh {\v^{\star}}
\def \vb {\v_{\star}}
\def \w {{w}}
\def \wh {\w^{\star}}
\def \wb {\w_{\star}}

\def \IS {\int_{\S}}
\def \IR {\int_{\R^3}}

\def\ep{\e}
\def\cn{\cdot \n}

\def \cL {\mathcal{L}}
\def \M {\mathbf{F}_1}

\def \dis {\displaystyle}

\title[Granular gases driven by inelastic scattering]
{\textbf{Equilibrium solution to the inelastic Boltzmann equation
driven by a particles thermal bath}}

\author{Marzia  Bisi, Jos\'{e}  A. Carrillo \& Bertrand Lods}

\address{\textbf{Marzia Bisi}, Dipartimento di Matematica, Universit\`{a} di Parma,
    Viale G. P. Usberti 53/A, 43100 Parma, Italia}    \email{marzia.bisi@unipr.it}

\address{\textbf{Jos\'{e} Antonio Carrillo}, ICREA (Instituci\'o Catalana de Recerca i Estudis Avan\c
cats) and Departament de Ma\-te\-m\`a\-ti\-ques, Universitat
Aut\`onoma de Barcelona, E-08193 Bellaterra, Spain}
\email{carrillo@mat.uab.es}

\address{\textbf{Bertrand Lods}, Laboratoire de Math\'{e}matiques, CNRS UMR 6620, Universit\'{e} Blaise Pascal (Clermont-Ferrand 2), 63177 Aubi\`{e}re Cedex,
France.}\email{bertrand.lods@math.univ-bpclermont.fr}

\begin{document}

\maketitle

\begin{abstract}
We show the existence of smooth stationary solutions for the
inelastic Boltzmann equation under the thermalization induced by a
host-medium with a fixed distribution. This is achieved by
controlling the $L^p$-norms, the moments and the regularity of the
solutions for the Cauchy problem together with arguments related
to a dynamical proof for the existence of stationary states.
\end{abstract}

\medskip

\section{Introduction}
\setcounter{equation}{0}

The dynamics of rapid granular flows is commonly modelled by a
suitable modification of the Boltzmann equation for inelastic
hard-spheres interacting through binary collisions \cite{BrPo,Vi}.
As well-known, in absence of energy supply, inelastic hard spheres
are cooling down and the energy continuously decreases in time. In
particular, the Boltzmann collision operator for inelastic hard
spheres does not exhibit any non trivial steady state. This is no
more the case if the spheres are forced to interact with an
external agency (thermostat) and, in such a case, the energy
supply may lead to a non-equilibrium steady state. For such driven
system (in a space homogeneous setting), the time-evolution of the
one-particle distribution function \(f(v,t)\), \(v\in\R^3\),
\(t>0\) satisfies the following
\begin{equation}
\label{be:force} \partial_t f  = \tau\Q(f,f) + \G(f),
\end{equation}
where $\tau \geq 0$ is a given constant, \(\Q(f,f)\) is the
inelastic Boltzmann collision operator, expressing the effect of
binary collisions of particles, while \(\G(f)\) models the forcing
term.

There exist in the literature several physical possible choices for
the forcing term $\G$ in order to avoid the cooling of the granular
gas: stochastic heating, particles heating or scaled variables to
study the cooling of granular systems and even a nonlinear forcing
term given by the quadratic elastic Boltzmann operator has been
taken into account~\cite{FoMi}. These options have been studied
first in the case of inelastic Maxwell models
\cite{Bobylev-Carrillo-Gamba,Carrillo-Cercignani-Gamba,Bobylev-Cercignani,
Bobylev-Cercignani2,SpTo,BisiCT,BisiCT2, CaTo,
Bobylev-Cercignani-Gamba2}. The most natural one is the pure
diffusion thermal bath for which
 \begin{equation}
\label{force:Gauss} \G (f) = \mu \,\Delta {f}
\end{equation}
where \(\mu>0\) is a constant, studied in \cite{GPV**,MiMo3} for
hard-spheres. Such a forcing term corresponds to the physical
situation in which granular beads receive random kicking in their
velocity, like air-levitated disks \cite{BSSS}. Another example is
the thermal bath with linear friction
\begin{equation}
\label{force:FP} \G(f)= \mu \,\Delta {f} + \lambda \mathrm{div}
(v\,f),
\end{equation}
where \(\lambda\) and \(\mu\) are positive constants. We also have
to mention the fundamental example of anti-drift forcing term which
is related to the existence of self-similar solution to the
inelastic Boltzmann equation:
\begin{equation}
\label{force:selfsimi} \G (f) = -\kappa \,\mathrm{div} (v f), \qquad
\kappa>0.
\end{equation}
This problem has been treated in \cite{MMR,MiMo,MiMo2} for
hard-spheres. For all the forcing terms given by
\eqref{force:Gauss}, \eqref{force:FP}, \eqref{force:selfsimi} it
is possible to prove the existence of a non-trivial stationary
state $F \geq 0$ such that
$$
\tau\Q(F,F)+\G(F)=0.
$$
Moreover, such a stationary state can be chosen to be smooth, i.e.
$F \in C^\infty(\R^3)$. Finally, even if the uniqueness (in
suitable class of functions) of such a stationary state is an open
problem, it can be shown for all these models that, in the weakly
elastic regime in which the restitution coefficient is close to
unity, the stationary state is unique. For an exhaustive survey of
the ``state of art'' on the mathematical results for the evolution
of granular media see~\cite{Vi}.

We are concerned here with a similar question when the forcing term
$\G$ is given by a \textit{\textbf{linear scattering operator}}.
This corresponds to a situation in which the system of inelastic
hard spheres is immersed into a so called \textit{\textbf{particles
thermal bath}}, i. e.  $\G$ is given by a linear Boltzmann collision
operator of the form:
$$\G(f)=\mathcal{B}[f,\M]$$
where $\M$ stands for the distribution function of the host fluid
and $\mathcal{B}[\cdot,\cdot]$ is a given collision operator for
(elastic or inelastic) hard-spheres.  The precise definition of $\G$
is given in Subsection~2.1.

This kinetic model has already been tackled for instance in
\cite{Bisi, BiSp} in order to derive closed macroscopic equations
for granular powders in a host medium. Let us also mention the
work \cite{biben} that investigates the case of a particles
thermal bath made of {\it elastic} hard-spheres at thermodynamical
equilibrium (i.e. $\M$ is a suitable Maxwellian).  The deviations
of the steady state (which is there assumed to exist) from the
Gaussian state are analyzed numerically. For inelastic Maxwellian
molecules, the existence of a steady state for a particles thermal
bath has been obtained in \cite{CaTo}.  To our knowledge, the
existence of a stationary solution of \eqref{be:force} for
particles bath heating and \textit{inelastic hard-spheres}  is an
open problem and it is the main aim of this paper.

Our strategy, inspired by several works in the kinetic theory of
granular gases \cite{GPV**,MiMo} or for coagulation-fragmentation
problems \cite{BaLa,Esco}, is based on a dynamic proof of the
existence of stationary states, see \cite[Lemma 7.3]{CaTo} for a
review. The exact ``fixed point theorem'' used here is reported in
Subsection 2.2. The identification of a suitable Banach space and
of a convex subset that remains invariant during the evolution, is
achieved by controlling moments and $L^p$--norms of the solutions.
In Section 3, we present regularity properties of the gain part of
both collision operators $\Q$ and $\G$ in~(\ref{be:force}). Then,
in Section 4 we get at first uniform bounds for the moments and
the Lebesgue norms; in addition, we prove the strong continuity of
the semi--group associated to~(\ref{be:force}), and the existence
and uniqueness of a solution to the Cauchy problem. All this
material allows to obtain, in Section~5, existence of non--trivial
stationary states. Finally, Section~6 contains the study of
regularity of stationary solutions. Many technical estimates
involving the quadratic dissipative operator $\Q(f,f)$ are based
on results presented in \cite{BouchutDesvillettes, MiMo,
MiMo2,MoVi} and in the references therein, but their extension to
the linear inelastic operator $\G(f)$ is not trivial at all for
the following reasons. First, since $\G$ is not quadratic, it
induces a lack of symmetry particularly  relevant in the study of
propagation of $L^p$ norms. Second, since the microscopic
collision mechanism is affected by the mass ratio of the two
involved media (thermal bath and granular material),
\textit{Povzner-like estimates} for $\G$ are not straightforward
consequences of previous results from \cite{GPV**}. Let us finally
mention that our analysis also applies to  linear scattering model
which corresponds to the case $\tau = 0$. For such a linear
Boltzmann operator, we obtain the existence of an equilibrium
solution, generalizing the results of \cite{MaPi,LoTo,SpTo} to
non-necessarily Maxwellian host distribution.

\section{Preliminaries}

Let us introduce the notations we shall use
in the sequel. Throughout the paper we shall use the notation
$\langle \cdot \rangle = \sqrt{1+|\cdot|^2}$. We denote, for any
$\eta \in \R$, the Banach space
\[
     L^1_\eta = \left\{f: \R^3 \to \R \hbox{ measurable} \, ; \; \;
     \| f \|_{L^1_\eta} := \int_{\R^3} | f (v) | \, \langle v \rangle^\eta \d\v
     < + \infty \right\}.
\]

More generally we define the weighted Lebesgue space $L^p_\eta
(\R^3)$ ($p \in [1,+\infty)$, $\eta \in \R$) by the norm
    \[ \| f \|_{L^p_\eta(\R^3)} = \left[ \int_{\R^3} |f (v)|^p \, \langle v
        \rangle^{p\eta} \d\v \right]^{1/p}. \]

The weighted Sobolev space $W^{k,p}_\eta (\R^3)$ ($p \in
[1,+\infty)$, $\eta \in \R$ and $k \in \mathbb{N}$) is defined by
the norm
    \[ \| f \|_{W^{k,p}_\eta (\R^3)}  =   \left[ \sum_{|s| \le k} \|\partial_v^s f\|_{L^p_\eta} ^p \right]^{1/p} \]
where $\partial_v^s$ denotes the partial derivative associated with
the multi-index $s \in \mathbb{N}^N$. In the particular case $p=2$
we denote $H^k _\eta=W^{k,2} _\eta$. Moreover this definition can be
extended to $H^s _\eta$ for any $s \ge 0$ by using the Fourier
transform.

\subsection{The kinetic model}

  We assume
the granular particles to be perfectly smooth hard spheres of mass
$m=1$ performing inelastic collisions. Recall that, as usual, the
inelasticity of the collision mechanism is characterized by a single
parameter, namely the coefficient of normal restitution \(0 <
\epsilon < 1\).  To define the collision operator we write
\begin{equation}
\label{co:gain-loss} \Q(f,f) = \Q^+(f,f) - \Q^-(f,f),
\end{equation}
where the ``loss'' term \(\Q^-(f,f)\) is
\begin{equation}
\label{co:loss} \Q^-(f,f) = f(f \ast |v|),
\end{equation}
and the ``gain'' term \(\Q^+(f,f)\) is given by
$$
\Q^+(f,f) = \frac{1}{4\pi\epsilon^2} \IR \int_{\mathbb{S}^2} |v-w|
 f('v) f('w) \,\d\sigma\d\w,
$$
where the pre-collisional velocities read as
\begin{equation}
\label{co:transfpre}
  'v=v+\frac{\zeta}{2\epsilon}\,(|q|\sigma-q),
\qquad
'w=w-\frac{\zeta}{2\epsilon}\,(|q|\sigma-q),
\end{equation}
with \(\zeta=\frac{1+\epsilon}{2}\). Notice that we always have
\(\frac{1}{2}<\zeta<1\). Its weak formulation will be the main
tool in the rest and it reads as
\begin{equation}
\label{co:gainw} \IR \Q^+ (f,f)(v)\, \psi(v)\d\v  =
\frac{1}{4\pi}\IR\IR f(v)\,f(w)\,|q|\, \IS \psi(v')\d\sigma\,
\d\w\d\v,
\end{equation}
where \(q=v-w\) is the relative velocity of two particles about to
collide, and \(v'\) is the velocity after the collision. The
collision transformation that puts \(v\) and \(w\) into
correspondence with the post-collisional velocities \(v'\) and
\(w'\) can be expressed as follows:
\begin{equation}
\label{co:transf}
  v'=v+\frac{\zeta}{2}\,(|q|\sigma-q),
\qquad  w'=w-\frac{\zeta}{2}\,(|q|\sigma-q).
\end{equation}

Combining \eqref{co:loss} and \eqref{co:gainw} and using the
symmetry that allows us to exchange \(v\) with \(w\) in the
integrals we obtain the following symmetrized weak form
\begin{equation}\label{co:weak}
\begin{split}
\IR \Q (f,f)(v)\, \psi(v)\d\v  =  \frac{1}{2} \IR\IR f(v)\,f(w)\,|q|
\mathcal{A}_\zeta[\psi](v,w)\d\w\d\v,
\end{split}
\end{equation}
where
\begin{equation}
\label{coll:psi} \mathcal{A}_\zeta[\psi](v,w) =
\frac{1}{4\pi}\IS(\psi(v')+\psi(w')-\psi(v)-\psi(w))\d{\sigma}.
\end{equation}

The inelastic Boltzmann operator $\Q (f,f)$ satisfies the basic
conservation laws of mass and momentum, obtained by taking
$\psi=1,v$ in the weak formulation \eqref{co:weak}, since
$\mathcal{A}_\zeta[1]=\mathcal{A}_\zeta[v]=0$. On the other hand, in
the modelling of dissipative kinetic equations, conservation of
energy does not hold. In fact, we obtain
$\mathcal{A}_\zeta[|v|^2]= -\frac{1-\epsilon^2}{4} |v-w|^2$ from
which we deduce
\begin{equation}
\IR \Q(f,f)(v)\,|v|^2 \d\v = -\frac{1-\epsilon^2}{8} \IR \IR |v-w|^3
f(v) f(w) \d\v\d\w , \label{disenergy}
\end{equation}
where we observe the dissipation of kinetic energy. In the absence
of any other source of energy, the system cools down as $t\to
\infty$ following Haff's law as proved in \cite{MiMo}.

As already said in Introduction, the forcing term $\G$ arising in
the kinetic equation~(\ref{be:force}) is chosen to be a linear
scattering operator, corresponding to the so called
\textit{\textbf{particles bath heating}},
\begin{equation}\label{linearoperator}
\G(f):=\cL (f)=\dfrac{1}{2\pi \lambda}\IR \IS  |q \cdot \n| \left[
 e^{-2}f( \vb )\M(\wb) -f( \v )\M(\w)\right]\d\w \d\n
\end{equation}
where $\lambda$ is the mean free path, $q=\v-\w$ is the relative
velocity, $\v_\star$ and $\w_\star$ are the pre-collisional
velocities which result, respectively, in $ \v $ and $\w$ after
collision. The collision mechanism related to the linear
scattering operator is characterized by
\begin{equation}\label{coef}
(\v-\w ) \cn = - \e(\vb-\wb)\cn ,
\end{equation}
where $\n \in \S$ is the unit vector in the direction of impact and
$0 < \e < 1$ is the \textit{constant restitution coefficient}
(possibly different from $\epsilon$). Here, we will consider a
similar separation of the operator into gain and loss terms, $\cL
(f)=\cL^+ (f)-\cL^- (f)$, with obvious definitions. Here the host
fluid is made of hard-spheres of mass $m_1$ (possibly different from
the traced particles mass $m=1$) and the distribution function $\M$
of the host fluid fulfils the following:

\begin{hyp}\label{AssF1} $\M$ is a nonnegative normalized distribution function  with bulk velocity $\u_1 \in \R^3$ and
temperature $\Theta_1>0$. Moreover, $\M$ is smooth in the following
sense,
$$\M \in H^s_\delta(\R^3), \qquad \forall s,\delta \geq 0$$
and of finite entropy
$\ds\int_{\R^3} \M(\v) \log \M(\v) \d v < \infty.$
\end{hyp}
\begin{nb}
Notice that, since $\M \in L^1_2$ is of finite entropy, it is
well-known \cite[Lemma 4]{arkinfty} that there exists some $\chi >0$
such that
\begin{equation}\label{hypp}
\nu(v):= \frac{1}{2\pi\, \lambda}\, \IR \IS |(\v-w) \cdot \n|\M(\w)\d
\w\d \n \geq \chi  \sqrt{1+|v|^2} \qquad \forall v \in \R^3.
\end{equation}
A particular choice of the distribution function $\M$, corresponding
to a host fluid at thermodynamical equilibrium, is the following
Maxwellian distribution
\begin{equation}\label{maxwe1}
\M(\v)=\mathcal{M}_1(v)=\bigg(\dfrac{m_1}{2\pi
\Theta_1}\bigg)^{3/2}\exp
\left\{-\dfrac{m_1(\v-\u_1)^2}{2\Theta_1}\right\}, \qquad \qquad \v
\in \R^3,
\end{equation}
Notice however that our approach remains valid for more general
distribution function.\end{nb}

For particles of mass $m=1$ colliding inelastically with
particles of mass $m_1$, the restitution coefficient being
constant, the expressions of the pre-collisional velocities
$(\v_\star,\w_\star)$ are given by~\cite{BrPo,SpTo}
\begin{equation*}
\v_{\star}=\v-2\alpha\dfrac{1-\beta}{1-2\beta} \left(\q \cdot
\n\right) \n,\qquad
\w_{\star}=\w+2(1-\alpha)\dfrac{1-\beta}{1-2\beta} \left(\q \cdot
\n\right)  \n,
\end{equation*}
where $\alpha$ is the mass ratio and $\beta$ denotes the
inelasticity parameter
\[
\alpha = \frac{m_1}{1+m_1} \in (0,1), \qquad \ \ \ \beta=\frac{1-\e}{2}
\in [0,1/2).
\]
The post-collisional velocities are given by
  \begin{equation}\label{postco}
  \v^\star = \v - 2\alpha (1-\beta)\left(\q \cdot
\n\right)
  \n,\qquad \w^\star = \w + 2(1-\alpha) (1-\beta)\left(\q \cdot
\n\right)  \n.
  \end{equation}
This linear operator can also be represented in a form closer to
\eqref{co:weak}. By making use of the following identity
\cite{Boby75,Desvillettes},
$$
\int_{\S} (\hat{q}\cdot \n)_+ \ \varphi (\n(q\cdot \n))\d\n = \frac14
\int_{\S} \varphi \left( {{q-|q|\sigma}\over2} \right)\d\sigma
$$
for any function $\varphi$, with $\hat{q}=q/|q|$, we can rewrite
the operator as
\begin{equation}\label{linearoperator2}
\cL (f)=\dfrac{1}{4\pi \lambda}\IR \IS  |q| \left[
 e^{-2}f(\tilde{v}_{\star})\M(\tilde{w}_{\star}) -f( \v )\M(\w)\right]\d\w \d\sigma
\end{equation}
with
\begin{equation*}
 \tilde{v}_{\star}=\v-\alpha\dfrac{1-\beta}{1-2\beta} \left(\q
-|\q|\sigma\right), \qquad
\tilde{w}_{\star}=\w+(1-\alpha)\dfrac{1-\beta}{1-2\beta}\left(\q
-|\q|\sigma\right).
\end{equation*}
For such a description, the post-collisional velocities are
\begin{equation}\label{postco2}
  \tilde{v}^\star = \v - \alpha (1-\beta)\left(\q
-|\q|\sigma\right),\qquad \tilde{w}^\star = \w + (1-\alpha)
(1-\beta)\left(\q -|\q|\sigma\right).
  \end{equation}
We consider Eq. \eqref{be:force} in the weak form: for any regular
$\psi=\psi(v)$, one has
\begin{multline}\label{weak}
\dfrac{\d }{\d t}\IR f(v,t)\psi(v)\d v=\frac{\tau}{2} \IR \IR
f(v,t)\,f(w,t)\,|q| \mathcal{A}_\zeta[\psi](v,w)\d\w\d\v \\
+\frac{1}{\lambda} \,
       \IR \IR |q|
       f(\v,t)\M(\w)\mathcal{J}_e[\psi](v,w) \d
       \v \d \w
\end{multline}
where
$$
\mathcal{J}_e[\psi](v,w)=\dfrac{1}{2\pi}\int_{\S}
|\hat{q}\cdot\n|\left(\psi(v^\star)-\psi(v)\right)\d\n=
\dfrac{1}{4\pi}\int_{\S}
\left(\psi(\tilde{v}^\star)-\psi(v)\right)\d\sigma.
$$

\subsection{Proof of stationary states: basic tools and strategy}

As stated in the Introduction, the final purpose of this paper is to prove the existence of a
non-trivial regular stationary solution $F \geq 0$ to~\eqref{be:force}. Namely, we look for $F \in L^1$, $F \geq 0$ such
that
\begin{equation}\label{stationaire}
\tau\Q(F,F)+\G(F)=0.
\end{equation}
\begin{nb} Notice that such a problem is trivial in the elastic case $\epsilon=1$ and whenever $\M$ is the Maxwellian
distribution \eqref{maxwe1}. Indeed, in such a case,  the Maxwellian
equilibrium distribution $\mathcal{M}^\sharp$ of $\cL$ provided by
\cite{MaPi,LoTo,SpTo} is a stationary solution to~(\ref{be:force})
since $\Q(\mathcal{M}^\sharp,\mathcal{M}^\sharp)=0$ (elastic
Botlzmann equation) and $\cL(\mathcal{M}^\sharp)=0.$\end{nb}

 The main ingredients are to show the existence of fixed
points for the flow map at any time, and thus a continuity in time
argument of the semi-group that allows to identify this
one-parameter family of fixed points as a stationary point of the
flow. Contraction estimates were used in \cite{BisiCT,CaTo} while in
the hard-sphere case the Tykhonov Fixed Point Theorem was the tool
needed \cite{GPV**,Esco,MiMo,BaLa}.

\

The exact result that will be used can be summarized as:

\begin{lemme}[\textit{\textbf{Dynamic proof of stationary states}}]\label{GPV}
Let $Y$ be a Banach space and $(S_t)_{t \ge 0}$ be a continuous
semi-group on $Y$ such that
\begin{enumerate}
\item[i)] there exists $Z$ a nonempty convex and weakly
(sequentially) compact subset of\, ${Y}$ which is invariant under
the action of $S_t$  (that is $S_t z  \in {Z}$ for any $z \in {Z}$
and $t \ge 0$);

\item[ii)] $S_t$ is weakly (sequentially) continuous on $Z$ for
any $t>0$.
\end{enumerate}
Then there exists $z_0 \in Z$ which is stationary under the action
of $S_t$ (that is $S_t z_0=z_0$ for any $t \ge 0$).
\end{lemme}

The strategy is therefore to identify a Banach space $Y$ and a
convex subset $Z \subset  Y$ in order to apply the above result. To
do so, one shall prove that
\begin{itemize}
\item for any $f_0 \in  {Y}$, there is a solution $f \in
\mathcal{C}\left([0,\infty), {Y}\right)$ to Eq. \eqref{be:force}
with $f(t=0)=f_0$; \item the solution $f$ is unique in $ {Y}$ and
if $f_0 \in  {Z}$ then $f(t) \in  {Z}$ for any $t \geq 0$;
\item the set $ {Z}$ is (weakly sequentially) compactly embedded
into $Y$;
\item solutions to \eqref{be:force} have to be (weakly sequentially)
stable, i.e.,  for any sequence $(f_n)_n \subset
\mathcal{C}\left([0,\infty), {Y}\right)$ of solutions to
\eqref{be:force} with $f_n(t) \in  {Z}$ for any $t \geq 0$, then,
there is a subsequence $(f_{n_k})_k$ which converges weakly to some
$f \in \mathcal{C}\left([0,\infty),  {Y}\right)$ such that $f$ is a
solution to \eqref{be:force}.
\end{itemize}

If all the above points are satisfied by the evolution problem
\eqref{be:force}, then one can apply Lemma \ref{GPV} to the
semi-group $(\mathcal{S}_t)_{t \geq 0}$ which to any $f_0 \in {Y}$
associates the unique solution $f(t)=\mathcal{S}_t f_0$
to~\eqref{be:force}. Moreover, the regularity properties of the gain
part of the operators \cite{MiMo} shall provide us the needed
regularity to show the existence of smooth stationary states.

\section{Regularity of gain operators}


We recall the following result, taken from \cite[Theorem 2.5,
Proposition 2.6]{MiMo} and based on \cite{BouchutDesvillettes,Lu},
on the regularity properties of the gain part operator $\Q^+(g,f)$
that we state here only for hard-spheres interactions in space dimension $N=3$.

\begin{propo}[\textit{\textbf{Regularity of the gain term $\Q^+$}}]\label{Q^+}
For all $s,\eta>0$, we have
\[
     \left\|\Q^+(g,f)\right\|_{H^{s+1}_\eta} \le
     C (s,\eta,\epsilon) \left[\left\|g\right\|_{H^{s}_{\eta+2}} \left\|f\right\|_{H^{s}_{\eta+2}} +
     \left\|g\right\|_{L^1_{\eta+2}} \left\|f\right\|_{L^1_{\eta+2}}\right]
\]
where the constant $C (s,\eta,\epsilon)>0$ only depends on the
restitution coefficient $\epsilon \in(0,1]$, $s$ and $\eta$.
Moreover, for any $p \in [1,\infty)$ and $\delta >0$,  there exist
$\theta\in (0,1)$ and a constant $C_\delta >0$, only depending on
$p$, $\epsilon$ and $\delta$, such that
       \[ \int_{\R^3} \Q^+(f,f) \, f^{p-1} \d\v \le
          C_\delta \,  \|f\|_{L^1} ^{1+p\theta} \, \|f\|_{L^p} ^{p(1-\theta)}
                         + \delta \, \|f\|_{L^1_2} \, \|f\|_{L^p _{1/p}} ^p .\]
\end{propo}

On the other hand, the linear operator $\cL(f)$ is quite similar
to the quadratic Boltzmann operator associated to hard-spheres
interactions and constant restitution coefficient $\ep$ by fixing
one of the distributions. In fact, it is possible to obtain the
following similar result:

\begin{propo}[\textit{\textbf{Regularity of the gain term $\cL^+$}}]\label{L^+}
For all $s,\eta>0$, we have
\begin{equation}\label{regl1}
     \left\|\cL^+(f)\right\|_{H^{s+1}_\eta} \le
     C (s,\eta,\ep) \left[\left\|\M\right\|_{H^{s}_{\eta+2}} \left\|f\right\|_{H^{s}_{\eta+2}} +
     \left\|\M\right\|_{L^1_{\eta+2}} \left\|f\right\|_{L^1_{\eta+2}}\right]
\end{equation}
where the constant $C (s,\eta,\ep)>0$ only depends on the
restitution coefficient $\ep \in(0,1]$, $s$ and $\eta$. Moreover,
 for any $p \in (1,\infty)$ and $\delta >0$,  there exist $q <
p$  and a constant $K_\delta >0$, only depending on $p$, $\e$ and
$\delta$, such that
\begin{multline}\label{regl2}
       \int_{\R^3} \cL^+(f) \, f^{p-1} \d\v \le
          K_{ \delta} \,  \|\M\|_{L^q}\, \|f\|_{L^p}^{p-1}\,
          \|f\|_{L^1}+ \delta \, \bigg(\|\M\|_{L^1_2} \, \|f\|_{L^p_{1/p}}^p
          \\\
                        + \|\M\|_{L^p_{1/p}} \, \|f\|_{L^1_2}\, \|f\|_{L^p_{1/p}}^{p-1} \bigg).
\end{multline}
\end{propo}

\begin{proof}
The proof of these two estimates relies on the same steps as in
Sections 2.2, 2.3 and 2.4 of~\cite{MiMo}, see also \cite{MoVi}. We
need just to have the same basic estimates as in their case. We
start with the proof of \eqref{regl1}. An expression of the Fourier
transform of $\cL^+$ can be obtained~as:
$$
{\mathcal F}\left[\cL^+(f)\right](\xi):=\int_{\R^3} \exp(-i \xi
\cdot v) \cL^+(f)(\v)\d\v =\dfrac{1}{4\pi \lambda} \int_{\S}
\widehat{G}(\xi_+,\xi_-)\d\sigma
$$
with $G(v,w)=|v-w|f(\v)\M(w)$, $\widehat{G}$ its Fourier transform with
respect to $(v,w)$ and
\begin{equation*}
\xi_+= (1-\alpha (1-\beta)) \xi + \alpha(1-\beta)|\xi|\sigma, \qquad
\xi_-= \alpha (1-\beta) \xi - \alpha(1-\beta)|\xi|\sigma.
\end{equation*}
With this expression at hand, it is immediate to generalize to
$\cL^+$ the regularity result in \cite[Theorem 2.5, Proposition
2.6]{MiMo} giving \eqref{regl1}.

\medskip

Now, let us prove the second result. We first notice that, as in
\cite{ArLo}, the gain operator $\cL^+$ admits an integral
representation. Actually, even if it is assumed in \cite{ArLo} that
$\M$ is given by the Maxwellian distribution \eqref{maxwe1}, a
careful reading of the calculations of \cite{ArLo} yields
\begin{equation} \cL^+f(\v)= \int_{\R^3}f(\w)k(\v,w)\d\w,
\label{L+Carl}
\end{equation}
where
\begin{equation*}\label{k(v,w)}
k(\v,\w)=\dfrac{1}{2\ep^2\gamma^2|\v-\w|}\int_{V_2 \cdot
(\w-\v)=0}\M\left(\v+V_2+\dfrac{1-2\g}{2\gamma}(\w-\v)\right)\d
V_2\end{equation*} with $\gamma=\alpha \frac{1-\b}{1-2\b}$ and
$\g=(1-\alpha)\frac{1-\b}{1-2\b}$. Arguing as in \cite{MiMo}, we
define  the operator $\mathcal{T}$ related to the Radon transform:
$$
\mathcal{T}\::\:g \in L^1(\R^3,\d\v)\mapsto
\mathcal{T}g(v)=\dfrac{1}{|v|}\int_{z \bot v}g(\mu v + z)\d z
$$
where  $\mu=1 - \frac{1-2\g}{2\gamma}$. For any $h \in \R^3$, let
$\tau_{h}$ denote the translation operator $\tau_h f (v) =
f(v-h),$ for any $v \in \R^3.$ Then, for any $g \in L^1(\R^3,\d\v)$,
one sees that
\begin{equation*}
\begin{split} (\tau_{\w} \circ
\mathcal{T})(g)(v)&=\dfrac{1}{|\v-\w|}\int_{z \bot (v-w)} g(\mu
(v-w)+z)\d z\\
&=\dfrac{1}{|\v-\w|}\int_{z \bot (v-w)}g\left(v-w +z
+\dfrac{1-2\g}{2\gamma}(\w-\v)\right)\d z, \qquad \forall v,w \in
\R^3.
\end{split}
\end{equation*}
Choosing $g=\tau_{-w} \M$ leads to
the following expression of the kernel $k(v,w)$:
$$
k(v,w)=\dfrac{1}{2\epsilon^2\gamma^2}\bigg[\tau_{w} \circ \mathcal{T}
\circ \tau_{-w}\bigg](\M)(v),\qquad v,w \in \R^3.
$$

This previous computation is at the heart of the arguments of
\cite[Theorem 2.2]{MiMo}, from which one gets a version of Lions'
Theorem \cite{lions} for a suitable regularized cut-off kernel with
collision frequency of the form
$B_{m,n}(|q|,\hat{q}\cdot\sigma)=
\Phi_{S_n}(|q|)\,b_{S_m}(\hat{q}\cdot\sigma)$,
with $\Phi_{S_n}$
smooth and with compact support $\left[ \frac{2}{n}, n \right]$,
and $b_{S_m}$ smooth and supported in $\left[-1 + \frac{2}{m}, 1 - \frac{2}{m} \right]$. More precisely, defining the
smoothed-out operator in angular and radial variables
$\cL^+_{S_{m,n}}$ as in \cite[Section 2.4]{MiMo}:
\begin{equation}\label{linearoperatorreg}
\cL^+_{S_{m,n}} (f)=\dfrac{1}{4\pi \lambda\, e^2}\IR \IS
B_{m,n}(|q|,\hat{q}\cdot\sigma) f( \vb )\M(\wb) \d\w \d\sigma
\end{equation}
then, for any $\eta \in \R^+$ and any $p > 1$,
there is
$C(p,\eta,m,n)>0$ depending only on $p$, $\eta$ and $(m,n)$, such that
\begin{equation} \label{SobL}
\|\cL^+_{S_{m,n}}(f)\|_{L^p_\eta} \leq
C(p,\eta,m,n)\|\M\|_{L^q_\eta}\,\|f\|_{L^1_{2|\eta|}}
\end{equation}
for some $q < p$ given by $q=\frac{5p}{3 +2p}$ if $p \in (1,6]$
while $q=\frac{p}{3}$ {if $p \in [6,+\infty)$} (see \cite[Corollary~2.4]{MiMo}). In particular, H\"older's inequality leads to
$$\int_{\R^3} \cL^+_{S_{m,n}}(f) f^{p-1}\d\v \leq \left(\int_{\R^3}
f^p \d\v\right)^{\frac{p-1}{p}} \|\cL^+_{S_{m,n}}(f)\|_{L^p} \leq
 C( m,n) \|f\|_{L^1} \|\M\|_{L^q} \|f\|_{L^p}^{p-1}
$$
for some explicit constant $C( m,n) >0$.

Similarly, one can define
the remainder part of $\cL^+$ which splits as
$$ \cL^+-\cL^+_{S_{m,n}}=:\cL^+_{R_{m,n}}=\cL^+_{RS_{m,n}}+\cL^+_{SR_{m,n}}+\cL^+_{RR_{m,n}}$$ with
$$\cL^+_{RS_{m,n}} (f)=\dfrac{1}{4\pi \lambda\, e^2}\IR \IS
\Phi_{R_n}(|q|)\,b_{S_m}(\hat{q}\cdot\sigma) f( \vb )\M(\wb) \d\w
\d\sigma,$$
$$\cL^+_{SR_{m,n}} (f)=\dfrac{1}{4\pi \lambda\, e^2}\IR \IS
\Phi_{S_n}(|q|)\,b_{R_m}(\hat{q}\cdot\sigma) f( \vb )\M(\wb) \d\w
\d\sigma,$$
$$\cL^+_{RR_{m,n}} (f)=\dfrac{1}{4\pi \lambda e^2}\IR \IS
\Phi_{R_n}(|q|)\,b_{R_m}(\hat{q}\cdot\sigma) f( \vb )\M(\wb) \d\w
\d\sigma,$$ where $\Phi_{R_n}(|q|)=|q|-\Phi_{S_n}(|q|)$ and
$b_{R_m}(\hat{q} \cdot \sigma)= 1  -
b_{S_m}(\hat{q}\cdot \sigma)$, $q \in \R^3, \sigma \in \S$.
H\"older's inequality provides
$$
\int_{\R^3} \cL^+_{R_{m,n}}(f) \, f^{p-1} \d\v  \leq \|f\|_{L^p_{1/p}}^{p-1} \|\cL^+_{R_{m,n}}(f)\|_{L^p_{-1/p'}}
$$
with $p'$ such that $\frac{1}{p} + \frac{1}{p'} =1$, hence we have to estimate $L^p_\eta$ norms of
$\cL^+_{SR_{m,n}}$, $\cL^+_{RS_{m,n}}$, $\cL^+_{RR_{m,n}}$ for $\eta=-1/p'$.

One can easily use \cite[Theorem 2.1]{MiMo} to prove that, for any $\eta \in \R$,
$$\|\cL^+_{SR_{m,n}}(f) +\cL^+_{RR_{m,n}}(f)\|_{L^p_\eta} \leq
\varepsilon(m)\bigg( \|\M\|_{L^1_{|1+\eta|+|\eta|}} \|f\|_{L^p_{1+\eta}}
+ \|f\|_{L^1_{|1+\eta|+|\eta|}}\|\M\|_{L^p_{1+\eta}}\bigg)$$ for
some explicit constant $\varepsilon(m)$ that, since the angular part of the collision kernel is such that $\lim_{m \to
\infty}\|b_{R,m}\|_{L^1(\S)}=0$, converges to $0$ as $m$ goes
to infinity.

It remains to estimate the norm of
$\cL^+_{RS_{m,n}}(f)$. We follow now the lines of \cite[Chapter 9,
p.~395]{MoPhD} (which differs slightly from \cite[Proposition
2.6]{MiMo} and is more adapted to the linear case). Precisely, we
split $f$ as $f=f_r+f_{r^c}=f(v)\chi_{\{|v| \leq r\}}+f(v) \chi_{\{|
v | > r\}}$ for some $r >0$.
Then, as in \cite[p. 395]{MoPhD}, there is some positive constant $C >0$ such that
$$\|\cL^+_{RS_{m,n}}(f_r)\|_{L^p_\eta} \leq C \dfrac{r}{n}
\|\M\|_{L^1_{|2+\eta|+|\eta|}}\|f\|_{L^p_{1+\eta}}$$ while
$$\|\cL^+_{RS_{m,n}}(f_{r^c})\|_{L^p_\eta} \leq C \frac{m^\lambda}{r}
\|f\|_{L^1_{|2+ \eta| + |\eta|}}\|\M\|_{L^p_{1+\eta}}$$ with $\lambda
>0$.

Gathering all the above estimates
 we get, for $\eta=-1/p'$,
 \begin{multline*} \int_{\R^3} \cL^+_{R_{m,n}}(f) \, f^{p-1} \d\v
  \leq C
\|f\|_{L^p_{1/p}}^{p-1}\left(\frac{r}{n}\|\M\|_{L^1_2}\|f\|_{L^p_{1/p}}
+ \frac{m^\lambda}{r}\|f\|_{L^1_2}\|\M\|_{L^p_{1/p}}\right)
 \\
+\, \varepsilon(m) \left(\|\M\|_{L^1_1} \, \|f\|_{L^p_{1/p}}^p +
                         \|\M\|_{L^p_{1/p}} \, \|f\|_{L^1_1}\,
                         \|f\|_{L^p_{1/p}}^{p-1}\right) \\
\leq  \left( C\, \frac{r}{n} + \varepsilon(m) \right)\|\M\|_{L^1_2} \, \|f\|_{L^p_{1/p}}^p +
       \left( C\, \frac{m^\lambda}{r} + \varepsilon(m) \right)  \|\M\|_{L^p_{1/p}} \, \|f\|_{L^1_2}\,
                         \|f\|_{L^p_{1/p}}^{p-1}.
\end{multline*}
The proof follows then  by choosing first $m$ large enough then $r$
large enough and subsequently $n$ big enough.

\end{proof}

\section{Regularity estimates for the Cauchy
problem}\label{sec:existence}

\subsection{Evolution of mean velocity and temperature}\label{moments}
Let $f(v,t)$ be a nonnegative solution to \eqref{be:force}. Define the mass density,
the bulk velocity
$$
\varrho(t)=\IR f(v,t)\d v, \qquad \mathbf{u}(t)=\dfrac{1}{\varrho(t)}\int_{\R^3} v
f(v,t)\d v
$$
and the temperature
$$
\Theta(t)=\dfrac{1}{3\varrho(t)}\IR |v-\mathbf{u}(t)|^2
f(v,t)\d v, \qquad \forall t \geq 0.
$$
Note that Eq. \eqref{weak} for $\psi=1$ leads to the mass
conservation identity $\dot{\varrho}(t)=0$ i.e.
$$
\varrho(t)=\varrho(0):=1.
$$
Now, Eq. \eqref{weak} for $\psi(v)=v$ yields
$$
\dot{\mathbf{u}} (t) =-\dfrac{\a(1-\b)}{\lambda}\IR \IR
|v-w|(v-w)f(v,t)\M(w)\d v \d w, \qquad \forall t \geq 0
$$
which illustrates the fact that the bulk velocity is not
conserved. To estimate the second order moment of $f$, let us
introduce the auxiliary function:
$$
F(t)=\IR \IR |v-w|^2 f(v,t)\M(w)\d\v\d\w.$$
Notice that
\begin{equation} \label{Fdef}
F(t)= \IR |\v-\u_1|^2 f(\v,t) \d\v + \dfrac{3}{m_1}\Theta_1 = 3\,
\Theta(t) + |\u (t) - \u_1|^2 + \dfrac{3}{m_1}\Theta_1\,.
\end{equation}
In particular, to obtain uniform in time bounds of the mean
velocity and the temperature, it is enough to provide uniform in
time estimates of $F(t)$. With the special choice
$\psi(v)=|v-\u_1|^2$ one has
$$\mathcal{A}_\zeta[\psi](v,w)=\dfrac{\zeta(1-\zeta)|q|}{4\pi}\int_{\mathbb{S}^2}(\sigma\cdot
q-|q|)\d\sigma=-{\zeta(1-\zeta)|q|^2}=-\dfrac{1-\epsilon^2}{4}|q|^2$$
while
\begin{equation*}\begin{split}
\mathcal{J}_e[\psi](v,w)&= 2\a^2(1-\b)^2 |q|^2-2\a (1-\b)\langle q,
v-\mathbf{u}_1\rangle\\
&=-2\kappa(1-\kappa )|q|^2-2\kappa \langle q, w-\mathbf{u}_1\rangle,
\qquad v,w \in \R^3
\end{split}\end{equation*}
with $\kappa=\a(1-\b)=\dfrac{\a}{2}(1+\e) \in (0,1)$ and $\langle
\cdot, \cdot \rangle$ denoting the scalar product. It is easy to see
that
\begin{multline}\label{temp}\dot{F}(t)=- \dfrac{(1-\epsilon^2)\tau}{8}\IR \IR
f(v,t)f(w,t)|q|^3\d\v\d\w\\
-\frac{ 2\kappa(1-\kappa)}{ \lambda}\IR \IR |v-w|^3 f(v,t)\M(w)\d \v
\d\w \\+ \frac{2\kappa}{\lambda}\IR \IR |q|\langle q,
\mathbf{u}_1-w\rangle f(v,t)\M(\w)\d\v\d\w.\end{multline} Now, since
$\ds\IR f(v,t)\d v=1$ for any $t \geq 0$, Jensen's inequality yields
$$\IR f(w,t)|q|^3\d\w \geq \left|v-\IR
wf(\w,t)\d\w\right|^3=|v-\mathbf{u}(t)|^3$$ and consequently
\begin{equation*}
\begin{split}
\IR \IR f(v,t)f(w,t)|q|^3\d\v\d\w \geq \IR
|v-\mathbf{u}(t)|^3f(v,t)\d\v  \\ \geq \left(\IR
|v-\mathbf{u}(t)|^2f(v,t)\d\v\right)^{3/2}
= \Big( 3\, \Theta(t) \Big)^{3/2}
\end{split}
\end{equation*}
where we used again Jensen's inequality.
In the same way,
$$\IR \IR |q|^3 f(v,t)\M(w)\d\v\d\w \geq \left(\IR \IR
|v-\w|^2f(v,t)\M(\w )\d\v\d\w\right)^{3/2}=F(t)^{3/2}.$$
Finally, the third integral in \eqref{temp} is estimated as
\begin{equation*}\begin{split}
\IR \IR |q| &\langle q,\u_1-w\rangle f(v,t)\M(w)\d\v\d\w \leq \IR \IR
|q|^2 |\u_1-\w|f(v,t)\M(w)\d\v\d\w\\
& \leq 2 \IR |v-\u_1|^2 f(v,t) \d\v\IR |w-\u_1|\M(w)\d\w + 2 \IR |w-\u_1|^3 \M(w)\d\w\\
&\leq C_0 F(t) \end{split}\end{equation*} where
$$
C_0 =2\, \max \left\{ \IR |\w-\u_1|\M (\w)\d\w\:,\; \dfrac{\ds\IR
|\w-\u_1|^3\M(\w) \d\w}{\ds\IR |\w-\u_1|^2\M(\w) \d \w} \right\}.
$$

In conclusion, we obtain
\begin{equation}\label{Ft}
\dot{F}(t) \leq -\dfrac{(1-\epsilon^2)\tau}{8} \big( 3 \Theta(t)
\big)^{3/2}-\dfrac{2\kappa(1-\kappa)}{
\lambda}F(t)^{3/2}+\dfrac{2C_0\kappa}{ \lambda}F(t)
 \leq -\gamma_1 F(t)^{3/2}+ \gamma_2 F(t)
 \end{equation}
where $\gamma_1=\dfrac{2\kappa(1-\kappa)}{\lambda}
>0$ and $\gamma_2=\dfrac{2C_0\kappa}{\lambda} >0$. A simple use of the
maximum principle shows that
$$F(t) \leq \max\left\{\left(\frac{\gamma_2 }{\gamma_1 }\right)^2, F(0)\right\}, \qquad
\forall t \geq 0.$$ Because of \eqref{Fdef}, this leads to explicit
upper bounds of the temperature $\Theta(t)$ and the velocity
$|\mathbf{u}(t)-\mathbf{u}_1|,$ namely
\begin{equation}\label{thetat}
 \sup_{t \geq 0} \bigg(3\Theta(t) +|\u(t)-\u_1|^2\bigg) \leq \max\left\{\left(\frac{\gamma_2 }{\gamma_1 }\right)^2, F(0)\right\} <\infty. \end{equation}

\subsection{Propagation of moments} To extend the previous basic
estimates, in the the spirit of \cite{BoGaPa}, we deduce from
Povzner-like estimates some useful inequalities on the moments
$$\mathds{Y}_r(t)=\int_{\R^3} f(v,t)|v|^{2r}\d\v,\qquad t \geq 0,\: r \geq 1$$
where $f(t)$ is a solution to \eqref{be:force} with unit mass. One
sees from \eqref{be:force} that
$$\dfrac{\d }{\d t}\mathds{Y}_r(t)=\tau
 Q_r(t)+ L _r(t),$$ where
$$ Q_r(t)=\int_{\R^3}\Q(f,f)(v,t)|\v|^{2r}\d\v, \qquad  {L}_r(t)=\int_{\R^3}\cL(f)(v,t)|\v|^{2r}\d\v.$$
The calculations provided in \cite{GPV**,BoGaPa} allow  to estimate,
in an almost optimal way, the quantity $ Q_r$. One has to do the
same for $L_r(t) $ given by
$$L_r(t)=\frac{1}{\lambda}\IR \IR f(v,t)\M(w)
|\v-\w|\mathcal{J}_e[\,|\cdot|^{2r}](v,w)\d\v\d\w.$$

 To do so, let us derive
\textit{Povzner-like estimates} for $\cL$ in the spirit of
\cite{GPV**}. The application of the result of \cite{GPV**} is not
straightforward since, obviously, $\cL$ is not quadratic and
because of the influence of the mass ratio
$\alpha=\frac{m_1}{m+m_1}$ in the collision mechanism. Here, we
will write the mass of particles $m$ even if taken as unity for
the sake of the reader. To be precise, we are looking for
estimates of
$$
\mathcal{J}_e[|\cdot|^{2r}](v,w)=\dfrac{1}{2\pi}\int_{\mathbb{S}^2}|\hat{q}\cdot
\n|\left(|v^\star|^{2r}-|v|^{2r}\right)\d\n, \qquad r \geq 1.
$$
To do so, it shall be convenient to write
\begin{equation}
\mathcal{J}_e[|\cdot|^{2r}](v,w)=\dfrac{1}{2\pi
m^r}\int_{\mathbb{S}^2} |\hat{q}\cdot \n|\left\{
\Psi\left(m|\v^\star|^2\right)-\Psi\left(m|\v|^2\right)\right\}\d\n
\label{Je2r}
\end{equation}
where $\Psi(x)=x^r$, $r \geq 1$. We adopt the strategy used in
\cite{GPV**} and write
\begin{equation}
\Psi\left(m|\v^\star|^2\right)-\Psi\left(m|\v|^2\right)=\mathfrak{q}_e(\Psi)(v,w)+\Psi(m_1|w|^2)-\Psi(m_1|\wh|^2)
\label{Psi}
\end{equation}
where
$$\mathfrak{q}_e(\Psi)(v,w)=\Psi\left(m|\v^\star|^2\right)+\Psi\left(m_1|\w^\star|^2\right)-\Psi\left(m|\v|^2\right)
-\Psi\left(m_1|\w|^2\right).$$ Now,
$$\mathfrak{q}_e(\Psi)(v,w)=\mathfrak{p}_e(\Psi)(\v,\w)-\mathfrak{n}_e(\Psi)(v,w)$$
with
$$
\begin{cases}
\mathfrak{p}_e(\Psi)(v,w)=\Psi\left(m|v|^2+m_1|w|^2\right)-\Psi\left(m|v|^2\right)-\Psi\left(m_1|w|^2\right)\\
\mathfrak{n}_e(\Psi)(v,w)=\Psi\left(m|v|^2+m_1|w|^2\right)-\Psi\left(m|\vh|^2\right)-\Psi\left(m_1|\wh|^2\right).
\end{cases}
$$
Applying \cite[Lemma 3.1]{GPV**} to the function $\Psi$ with
$x=m|v|^2$ and $y=m_1|w|^2$, we see that there exists $A >0$ such
that
\begin{equation}\label{frakpe}
\mathfrak{p}_e(\Psi)(v,w) \leq A \Big( m|\v|^2
\Psi'\left(m_1|w|^2\right)+m_1|w|^2 \Psi'\left(m|v|^2\right)\Big)
\end{equation}
while, since $\Psi$ is nondecreasing and $m|v|^2+m_1|w|^2 \geq
m|\vh|^2+m_1|\wh|^2$, there exists $b >0$ such that
\begin{equation*} \label{frakne}
\mathfrak{n}_e(\Psi)(v,w) \geq b\, m\,m_1 |\vh|^2\,|\wh|^2
\Psi''\left(m|\vh|^2+m_1|\wh|^2\right).
\end{equation*}
One can then write
$$\mathfrak{n}_e(\Psi)(v,w) \geq b\, \Delta(\vh,\wh) \left(m|\vh|^2+m_1|\wh|^2\right)^2
\Psi''\left(m|\vh|^2+m_1|\wh|^2\right)$$ where
$$\Delta(\vh,\wh)=\dfrac{m|\vh|^2\,\,m_1|\wh|^2}{\left( m|\vh|^2+m_1|\wh|^2 \right)^2}\,.$$
To estimate better the above term $\Delta(\vh,\wh)$, it will be
convenient to parametrize the post-collisional velocities in the
\textit{center of mass--relative velocity} variables, which, with
respect to the usual transformation (see e.g. \cite[Eq.
(3.10)]{GPV**}) depend on the masses $m$ and $m_1$. Namely, let us
set
$$\vh=\dfrac{z + m_1\ell |q|\varpi}{m+m_1}, \qquad \wh=\dfrac{z -
m \ell |q|\varpi}{m+m_1}$$ where $z=mv+m_1w$, $q=v-w$ and $\varpi$
is a parameter vector on the sphere $\mathbb{S}^2.$ The parameter~$\ell$ is positive and such that $\vh-\wh=\ell|\v-\w|\varpi.$ In
particular, one sees from the representation~\eqref{postco} that $0
< \ell \leq 1.$ In this representation, one has
$$|\vh|^2=\dfrac{1}{(m+m_1)^2}\bigg( |z|^2+m_1^2 \ell^2|q|^2+2\ell m_1|q||z|\cos
\mu\bigg)$$ and
$$|\wh|^2=\dfrac{1}{(m+m_1)^2}\bigg(|z|^2+m^2 \ell^2|q|^2-2\ell
m|q||z|\cos \mu\bigg),$$ where $\mu$ is the angle between $z$ and
$\varpi$. One has then
\begin{equation}\label{energy}
m|\vh|^2+m_1|\wh|^2=\dfrac{1}{m+m_1}\bigg(|z|^2+\ell^2m\,m_1|q|^2\bigg).\end{equation}
One can check that
\begin{equation*}\begin{split}
(m|\vh|^2)&\,(m_1|\wh|^2) = \frac{m\, m_1}{(m+m_1)^4} \bigg\{ \Big[
|z|^2 + \ell^2 m m_1 |q|^2 \Big]^2 - \Big[ |z|^2 - \ell^2 m m_1
|q|^2 \Big]^2 \cos^2 \mu\\&+ \Big[ \ell (m_1-m) |z| |q| + \Big(
|z|^2 - \ell^2 m m_1 |q|^2 \Big) \cos \mu \Big]^2
  - 4 \ell^2 m m_1 |z|^2 |q|^2 \cos^2 \mu \bigg\},\end{split}\end{equation*}
 i.e.
$$
(m|\vh|^2)\,(m_1|\wh|^2) \geq \frac{m\, m_1}{(m+m_1)^4} \Big[ |z|^2
+ \ell^2 m m_1 |q|^2 \Big]^2 (1- \cos^2 \mu)\,.$$ Therefore
$$
\Delta(\vh,\wh) \geq \frac{m\, m_1}{(m+m_1)^2}\, \sin^2 \mu\,.
$$
 We obtain then an estimate similar to the one obtained in~\cite{GPV**}. Moreover, it is easy to see from~\eqref{energy}
 that
$$m|\vh|^2+m_1|\wh|^2 \geq \ell^2 \left(m|v|^2+m_1|w|^2\right)$$
and, arguing as in \cite{GPV**}, there exists some
constant $\eta
>0$ such
 that
\begin{equation}\label{estimateNE}\mathfrak{n}_e[\Psi](v,w) \geq \eta \sin^2\mu
 \left(m|v|^2+m_1|w|^2\right)^2
 \Psi''\left(m|v|^2+m_1|w|^2\right).\end{equation}
This allows to prove the following:

\begin{lemme}[\textit{\textbf{Povzner-like estimates for $\cL$}}]\label{povzner} Let $\Psi(x)=x^r,$ $r >1$. Then, there exist positive
constants $k_r  $ and $A_r $ such that
$$|\v-\w|\mathcal{J}_e[|\cdot|^{2r}](v,w) \leq A_r
\bigg(|v||w|^{2r}+|v|^{2r}|w|\bigg)+ \frac{m_1^r}{m^r}|v-w||w|^{2r}-k_r
\bigg(|v|^{2r+1}+|w|^{2r+1}\bigg) ,$$ for any $v,w \in \R^3.$
\end{lemme}
\begin{proof}
Bearing in mind that $\mathcal{J}_e[|\cdot|^{2r}](v,w)$ is provided by~(\ref{Je2r}) and (\ref{Psi}),
first of all, since $\Psi\left(m_1|\wh|^2\right) \geq 0$, we note that
$$\Psi(m|\vh|^2)-\Psi(m|v|^2) \leq
\mathfrak{q}_e[\Psi](v,w)+\Psi\left(m_1|\w|^2\right)=\mathfrak{q}_e[\Psi](v,w)+m_1^r
|\w|^{2r}.$$
Then, integrating \eqref{frakpe} and \eqref{estimateNE} with respect
 to the angle $\n \in \mathbb{S}^2$, one obtains, as in \cite[Lemma 3.3.]{GPV**} and \cite[Lemma 3.4]{GPV**}, that there are
 ${A}_r$ and $k_r >0$ such that, for any $v,w \in \R^3$:
$$|\v-\w|\dfrac{1}{2\pi m^r}\int_{\mathbb{S}^2}\mathfrak{q}_e(\Psi)(v,w)|\hat{q}\cdot
\n|\d\n \leq {A}_r\bigg(|v||w|^{2r}+|v|^{2r}|w|\bigg) -k_r
\bigg(|v|^{2r+1}+|w|^{2r+1}\bigg),$$
and this concludes the proof.
\end{proof}

The above Lemma (restoring $m=1$) together with the known estimates for $Q_r(t)$
allow to formulate the following

\begin{propo}[\textit{\textbf{Propagation of moments}}]\label{propo:moments} Let $f(t)$ be a solution to \eqref{be:force} with unit mass. For any $r \geq
1$, let
$$\mathds{Y}_r(t)=\int_{\R^3} f(v,t)|v|^{2r}\d\v,\qquad t \geq 0.$$
Then, there are positive constants $\mathbf{A}_r,$ $\mathbf{K}_r$
and $\mathbf{C}_r$ that depend only on $r,\a,\b, \tau, \lambda$ and
the moments of $\M$ such that
$$\dfrac{\d }{\d t}\mathds{Y}_r(t) \leq \mathbf{C}_r
 +\mathbf{A}_r
 \mathds{Y}_{r}(t)-\mathbf{K}_r\mathds{Y}_r^{1+1/2r}(t),\qquad \qquad \forall t \geq 0.$$
As a consequence, if $\mathds{Y}_r(0) < \infty,$ then $\sup_{t \geq
0} \mathds{Y}_r(t) < \infty.$
\end{propo}
\begin{proof} Recall that
$\frac{\d }{\d t}\mathds{Y}_r(t)=\tau   Q_r(t)+ {L}_r(t),$ where
$$ Q_r(t)=\int_{\R^3}\Q(f,f)(v,t)|\v|^{2r}\d\v, \qquad {L}_r(t)=\int_{\R^3}\cL(f)(v,t)|\v|^{2r}\d\v.$$
According to \cite[Lemma 3.4]{GPV**}, there exist $\widetilde{A}_r
>0$ and $\widetilde{k}_r >0$ such that
$$Q_r(t) \leq \widetilde{A}_r
\mathds{Y}_{1/2}(t)\mathds{Y}_{r }(t) -\widetilde{k}_r
\mathds{Y}_{r+1/2}(t), \qquad t \geq 0.$$ Now, from Lemma \ref{povzner}
\begin{multline*}
\lambda L_r(t) \leq A_r M_{ r} \mathds{Y}_{1/2}(t) + A_r
M_{1/2}\mathds{Y}_r(t) + m_1^r \IR \IR
|v-w||w|^{2r}f(v,t)\M(w)\d\w\d\v \\
-k_r \mathds{Y}_{ r +1/2}(t)-k_r M_{r+1/2},
\end{multline*}
where $M_s=\ds\IR |\w|^{2s} \M(\w)\d\w,$ $s \geq 1.$ One has
$$\IR \IR
|v-w||w|^{2r}f(v,t)\M(w)\d\w\d\v \leq M_r
\mathds{Y}_{1/2}(t)+M_{r+1/2}$$ and, denoting $c_{1/2}:=\sup_{t \geq
0}\mathds{Y}_{1/2}(t) < \infty$, one has $$L_r(t) \leq \mathbf{C}_r
+\frac{A_r M_{1/2}}{\lambda} \mathds{Y}_{r}(t) -
\frac{k_r}{\lambda}\, \mathds{Y}_{r+1/2}(t)$$ where
$\mathbf{C}_r=\left(c_{1/2}A_rM_r + c_{1/2}m_1^r M_r + m_1^r
M_{r+1/2}\right)/\lambda$ is a positive constant depending only on
$\a,\b,\lambda,$ $r \geq 1$ and the moments of $\M$. Gathering all
these estimates leads to
$$\dfrac{\d}{\d t} \mathds{Y}_{r}(t) \leq \mathbf{C}_r
 +\mathbf{A}_r
 \mathds{Y}_{r}(t)-\mathbf{K}_r
\mathds{Y}_{r+1/2}(t)$$ where $\mathbf{A}_r=\tau \widetilde{A}_r
c_{1/2} + \frac{1}{\lambda} A_r M_{1/2} >0$ and $\mathbf{K}_r=\tau
\widetilde{k}_r +\frac{k_r}{\lambda} >0$.
 Now,  thanks
to the mass conservation and H\"{o}lder's inequality, one gets
$\mathds{Y}_{r+1/2}(t) \geq  \mathds{Y}_r^{1+1/2r}(t)$ which leads
to the desired result.
\end{proof}
\begin{nb}
We see from the definition of the positive constants $\mathbf{A}_r,$
$\mathbf{C}_r$ and $\mathbf{K}_r$ that the above Proposition still
holds true whenever $\tau=0 $ (i.e. for the linear problem).
\end{nb}

\subsection{Propagation of Lebesgue norms} Let us consider now an
initial condition $f_0 \in L^1_2 \cap L^p$ for some $1 < p <\infty.$
We compute the time derivative of the $L^p$ norm of the solution
$f(v,t)$ to \eqref{be:force}:
\begin{multline*}
\dfrac{1}{p}\dfrac{\d }{\d t} \int_{\R^3}
f^p(v,t)\d\v=\tau\int_{\R^3} \Q^+(f,f)f^{p-1}\d\v - \tau\int_{\R^3}
f^{p-1}\Q^-(f,f)\d\v \\
+\int_{\R^3}\cL^+(f)f^{p-1}\d\v-\int_{\R^3}\cL^-(f)f^{p-1}\d\v.\end{multline*}
Using the fact that $\ds \int_{\R^3} f^{p-1}\Q^-(f,f)\d\v \geq 0$
and $\cL^-(f)(v)=\nu(v)f(v)$ where the collision frequency $\nu(v)$
is given by
$$\nu(v)=\frac{1}{2\pi \lambda} \IR \IS |(\v-w) \cdot \n|\M(\w)\d \w\d
\n,$$
we obtain the estimate:
$$ \dfrac{1}{p}\dfrac{\d }{\d t} \int_{\R^3}
f^p(v,t)\d\v \leq \tau \int_{\R^3} \Q^+(f,f) f^{p-1}\d\v +
\int_{\R^3} \cL^+(f)f^{p-1}\d\v - \int_{\R^3} \nu (v)f^p(v,t)\d\v.$$
Using the  lower bound \eqref{hypp}, we get
\begin{equation*}\begin{split}
\dfrac{1}{p}\dfrac{\d }{\d t}\|f(t)\|_{L^p}^p
&\leq \tau \int_{\R^3} \Q^+(f,f) f^{p-1}\d\v + \int_{\R^3}
\cL^+(f)f^{p-1}\d\v - \chi \|f\|_{L^p_{1/p}}^p\,.
\end{split}\end{equation*}
Proposition \ref{Q^+} and the conservation of mass imply that, for any $\delta
>0$, there is $\theta >0$ and some $C_\delta$ such that
$$
\int_{\R^3} \Q^+(f,f) f^{p-1}(\v,t)\d\v \leq C_\delta
\|f(t)\|^{p(1-\theta)}_{L^p}+\delta
\|f(t)\|_{L^1_2}\,\|f(t)\|^{p}_{L^p_{1/p}}.
$$
Moreover, Proposition \ref{L^+} implies that, for any $\delta>0$,
$$
\int_{\R^3} \cL^+(f,f) \, f^{p-1} \d\v \le \,
          C_1 \|f(t)\|_{L^p}^{p-1}\,
          + C_2 \delta \, \left(\|f(t)\|_{L^p_{1/p}}^p +
                        \|f(t)\|_{L^1_2}\, \|f(t)\|_{L^p_{1/p}}^{p-1}\right),
$$
for some constants $C_1,C_2 >0$ that depend only on
$p,\delta, \eta, \alpha, \ep$ and the norms of $\M$ in the spaces
involved in \eqref{regl2}. Recall that there is some $M_2$ such that
$$\sup_{t \geq 0}\|f(t)\|_{L^1_2}=1 + \sup_{t \geq
0}\int_{\R^3}|v|^2f(v,t)\d\v \leq M_2 <\infty.
$$
Now, using Young's inequality, $xy^{p-1}\leq \frac1p x^p +
\frac{p-1}p y^p$, for any $x,y\geq 0$, we have
$$
\int_{\R^3} \cL^+(f,f) \, f^{p-1} \d\v \le \,
         C_1  \|f(t)\|_{L^p}^{p-1}\,
          + C_3 \delta \, \left(\|f(t)\|_{L^p_{1/p}}^p +
                        M_2^p
                          \right)$$
                          for some constant $C_3 >0$.
Collecting all the bounds above, we get the estimate
\begin{align*}
\dfrac{1}{p}\dfrac{\d }{\d t}\|f(t)\|_{L^p}^p \leq & \tau C_\delta
\|f(t)\|^{p(1-\theta)}_{L^p}+\delta \tau M_2
\|f(t)\|_{L^p_{1/p}}^p + C_1\,  \|f(t)\|_{L^p}^{p-1}\,\\
         & + C_3 \delta \, \left(\|f(t)\|_{L^p_{1/p}}^p +
                        M_2^p\right) -\chi
\|f(t)\|^p_{L^p_{1/p}}\\
\leq & \tau C_\delta \|f(t)\|^{p(1-\theta)}_{L^p} +
\left( \delta (\tau M_2 + C_3) - \frac{\chi}{2} \right)
\|f(t)\|_{L^p_{1/p}}^p \\
         & + C_1\,  \|f(t)\|_{L^p}^{p-1} + C_3 \delta \, M_2^p -\frac{\chi}{2} \|f(t)\|^p_{L^p} ,
\end{align*}
since $\|\cdot\|_{L^p_{1/p}} \geq \|\cdot\|_{L^p}$. Choosing now
$\bar{\delta}$ such that $\bar{\delta} (\tau M_2 +C_3) < \chi/2$, we get the
existence of positive constants $C_4$, $C_5$ and $C_6$ such that
$$
\dfrac{1}{p}\dfrac{\d }{\d t}\|f(t)\|_{L^p}^p \leq
C_4\|f(t)\|^{p(1-\theta)}_{L^p}+ C_5 \|f(t)\|_{L^p}^{p-1} -C_6
\|f(t)\|^p_{L^p} + C_3M_2^p\bar{\delta}.
$$
It is not difficult to get then that $ \sup_{t \geq 0}
\|f(t)\|_{L^p} < \infty.$ This can be summarized in the following

\begin{propo}[\textit{\textbf{Propagation of $L^p$-norms}}]\label{propaLP}
Let $p\in(1,\infty)$ and $f_0 \in L^1_2 \cap L^p$ with unit mass.
Then, the solution $f(t)$ to~\eqref{be:force} satisfies the
following uniform bounds
$$\sup_{t \geq 0}\left( \|f(t)\|_{L^1_2}+\|f(t)\|_{L^p}\right) < \infty.$$
\end{propo}

\begin{nb} Notice that the fact that $\M$ is of finite entropy (see Assumption
\ref{AssF1}) has been used here above, via the lower bound
\eqref{hypp},   in order to control from below $L^p$ norms involving
the loss operator $\cL^-$. Whenever $\tau >0$, it is possible then
to replace such estimates involving $\cL^-$ by others that
involve~$\Q^-$. Notice also that, whenever $\tau=0$ (i.e. in the
linear case), only the above constant $C_4$ vanishes and we still
have $\sup_{t \geq 0}\|f(t)\|_{L^p} < \infty$.
\end{nb} As a corollary, we deduce as in \cite[Section 3.4]{MiMo}, see also
\cite{CDT}, the following non-concentration result:

\begin{propo}[\textit{\textbf{Uniform non-concentration}}]
Let $f_0$ be given with unit mass. Assume that there exists some
$p\in(1,\infty)$ such that $f_0 \in L^1_2 \cap L^p$. Then, there
exists some positive constant $\nu_0$ such that
$$
 \nu_0 \leq \int_{\R^3} |\v-\u(t)|^2 f(v,t)\d v \leq
1/\nu_0, \qquad \forall t \geq 0,
$$
where $f(v,t)$ is the solution to \eqref{be:force} with $f(0)=f_0$
and $\u(t)=\ds\int_{\R^3} vf(v,t)\d\v,$ $t \geq 0.$
\end{propo}

\begin{proof} Let $f(t)$ be the solution to \eqref{be:force} with
$f(0)=f_0$. From the above Proposition, there exists $C_p >0$ such
that $\sup_{t \geq 0} \|f(t)\|_{L^p} \leq C_p,$ and H\"{o}lder's
inequality implies that, for any $r >0$,
$$\sup_{ t \geq 0} \int_{\{|\v-\u(t)| < r\}}f(v,t)\d\v \leq C_p\,
\left(\frac{4\pi}{3}r^3\right)^{\frac{p-1}{p}}.$$ Accordingly, there
is some $r_0
>0$ such that
$$\int_{\{|\v-\u(t)|< r_0\}} f(v,t)\d\v \leq \frac{1}{2}\,, \qquad \forall t
\geq 0.$$ Then, for any $t \geq 0$, recalling that $\dis \int_{\R^3}
f(v,t)\d\v=1$ for any $t \geq 0$,
\begin{equation*}\begin{split}
\int_{\R^3}f(v,t)|v-\u(t)| ^2\d\v &\geq \int_{\{|\v-\u(t)| \geq
r_0\}} f(v,t)| v-\u(t)|^2\d\v \geq
r_0^2 \int_{\{|\v-\u(t)| \geq r_0\}} f(v,t)\d\v\\
&\geq r_0^2\left(1-\int_{\{|\v-\u(t)| <
r_0\}}f(v,t)\d\v\right)\geq\frac{r_0^2}{2}
\end{split}\end{equation*}
which concludes the proof.
\end{proof}

\subsection{$L^1$-- stability}  As in \cite{MiMo}, in order to prove the strong
continuity of the semi-group $(\mathcal{S}_t)_{t \geq 0}$ associated
to \eqref{be:force}, one has to provide an estimate of
$\|f(t)-g(t)\|$ for two solutions $f(t)$ and $g(t)$ of
\eqref{be:force} with initial conditions $f(0), g(0)$ in some
subspace of $L^1$. This is the object of the following stability
result, inspired by \cite[Proposition 3.2]{MiMo} and
\cite[Proposition 3.4]{MMR}.

\begin{propo}[\textit{\textbf{$L^1$-stability}}]\label{stability} Let $f_0,g_0$ be two nonnegative functions of $L^1_3$
and let $f(t),g(t) \in \mathcal{C}(\R^+,L^1_2) \cap
L^\infty(\R^+,L^1_3)$ be the associated solutions to
\eqref{be:force}. Then, there is $\Lambda >0$ depending only on
$\sup_{t \geq 0} \|f(t)+g(t)\|_{L^1_3}$ such that
$$\|f(t)-g(t)\|_{L^1_2} \leq \|f_0-g_0\|_{L^1_2}\exp(\Lambda t),
\qquad \forall t \geq 0.$$
\end{propo}
\begin{proof} Let $h(t)=f(t)-g(t)$. Then, $h$ satisfies the following equation:
\begin{equation}\label{be:h}
\partial_t h(v,t)=\tau\bigg\{\Q(f,f)-\Q(g,g)\bigg\}+ \cL(h),\qquad
h(0)=f_0-g_0.\end{equation} As in \cite{MMR,MiMo}, the proof
consists in multiplying \eqref{be:h} by $\psi(
\v,t)=\mathrm{sgn}(h(v,t))\langle v \rangle ^2$ and integrating over
$\R^3$. We get
$$\dfrac{\d}{\d t}\int_{\R^3}|h(v,t)|\langle v\rangle^2\d\v =I(t)+L(t)$$
where
$$I(t)=\tau\int_{\R^3} \bigg\{\Q(f,f)-\Q(g,g)\bigg\}\psi(v,t)\d\v
\quad \text{  and } \quad
L(t)=\int_{\R^3}\cL(h)(v,t)\psi(v,t)\d\v.$$ To estimate the integral
$I(t)$ we resume the arguments of \cite[Proposition 3.4]{MMR} that
we shall need again later. According to \eqref{co:weak}
$$I(t) =\dfrac{\tau}{2}   \IR\IR \left(f(v,t)\,f(w,t)-g(v,t)g(w,t)\right)\,|q|
\mathcal{A}_\zeta[\psi(t)](v,w)\, \d\w\,\d\v. $$
The change of variables $(v,w) \mapsto (w,v)$ implies that
\begin{equation*}
I(t) =\frac{\tau}{2}\IR\IR
\left(f(v,t)-g(v,t)\right)\left(f(w,t)+g(w,t)\right) \,|q|
\mathcal{A}_\zeta[\psi(t)](v,w)\, \d\w\,\d\v.
\end{equation*}
Moreover, it is easily seen from the definition of $\psi$ that
\begin{equation*}\begin{split}
\left(f(v,t)-g(v,t)\right)&\mathcal{A}_\zeta[\psi(t)](v,w) \leq
\dfrac{1}{4\pi}\left| f(v,t)-g(v,t)\right|\IS \left(\langle v'
\rangle^2+\langle w'\rangle^2-\langle v\rangle^2+\langle w
\rangle^2\right)\,\d{\sigma}\\
&\leq 2\left|f(v,t)-g(v,t)\right|\langle w \rangle^2
\end{split}\end{equation*}
where we used the fact that $\mathcal{A}_\zeta[\langle \cdot
\rangle^2](v,w)=-\frac{1-\epsilon^2}{4}|q|^2 \leq 0.$ Therefore,
\begin{multline*} I(t) \leq  {\tau} \IR\IR  \,|q|\left|
f(v,t)-g(v,t)\right|\left(f(w,t)+g(w,t)\right) \langle w
\rangle^2\d\v\d\w\\
\leq \tau \IR \left|f(v,t)-g(v,t)\right| \langle v \rangle^2 \d\v
\IR \left(f(w,t)+g(w,t)\right) \langle w \rangle^3
\d\w\end{multline*} i.e. \begin{equation}\label{It} I(t) \leq \tau
\big \|f(t)+g(t)\big\|_{L^1_3} \big \|f(t)-g(t)\big \|_{L^1_2},
\qquad \qquad \forall t \geq 0.\end{equation} On the other hand,
recalling that $\cL(h)(v,t)=\cL^+(h)(v,t)-\nu(v)h(v,t)$, from formula~(\ref{L+Carl}) one has
\begin{equation*}\begin{split}
L(t)&=\int_{\R^3}\psi(v,t)\d\v\int_{\R^3}h(w,t)k(v,w)\d\w-\int_{\R^3}\nu(\v)|h(v,t)|\langle
v \rangle ^2 \d\v\\
&\leq\int_{\R^3}\langle v \rangle^2\d\v\int_{\R^3}|h(w,t)|k(v,w)\d\w
-\int_{\R^3}\nu(\v)|h(v,t)|\langle v \rangle ^2 \d\v,
\end{split}\end{equation*}
i.e. $L(t) \leq \ds\int_{\R^3}\cL(|h|)(v,t) \langle v \rangle ^2
\d\v.$ Now, since $\ds \int_{\R^3}\cL(|h|)(v,t)\d\v=0$ for any $h$,
one gets that
$$L(t) \leq \int_{\R^3}\cL(|h|)(v,t) \left|\v\right|^2 \d\v.$$
Resuming the calculations performed in Section \ref{moments} (see
Eq. \eqref{temp}), one gets that
\begin{multline*}\int_{\R^3}\cL(|h|)(v,t) \left|\v\right|^2 \d\v\leq
-\frac{ 2\kappa(1-\kappa)}{ \lambda}\IR \IR |v-w|^3 |h (v,t)|\M(w)\d \v
\d\w \\+ \frac{2\kappa}{\lambda}\IR \IR |q|\langle q, -w\rangle
|h(v,t)|\M(\w)\d\v\d\w\\
\leq \frac{2\kappa}{\lambda}\IR \IR
|q|^2|w||h(v,t)|\M(w) \d\v\d\w.\end{multline*} This leads to
$$L(t) \leq \frac{2\kappa}{\lambda}\bigg(2\IR  |v|^2 |h(v,t)|\d\v \IR |w|\M(w) \d\w + 2\IR |h(v,t)|\d\v \IR |w|^3\M(w)\d\w\bigg)$$
and, setting $c_+=\ds\frac{4\kappa}{\lambda} \max\left\{\IR
|w|\M(w)\d\w,\IR |w|^3\M(w)\d\w \right\},$  we get
$$L(t) \leq c_+ \int_{\R^3} |h(v,t)|\langle v \rangle^2 \d\v=c_+ \big\|f(t)-g(t)\big\|_{L^1_2}.$$
Gathering \eqref{It} together with the latter estimate and denoting
then $\Lambda= \tau\sup_{t \geq 0}\|f(t)+g(t)\|_{L^1_3} + c_+$, we
get the estimate
$$\dfrac{\d }{\d t}\big\|f(t)-g(t)\big\|_{L^1_2} \leq \Lambda \big\|f(t)-g(t)\big\|_{L^1_2}, \qquad \qquad t \geq 0$$
and the proof is achieved.
\end{proof}

\subsection{Well-posedness of the Cauchy problem}

We are in position to prove that the Boltzmann equation
\eqref{be:force} admits a unique regular solution in the following
sense:
\begin{theo}[\textit{\textbf{Existence and uniqueness of solution to the Cauchy problem}}] \label{CauchyPb}
Take an initial datum $f_0 \in L_3^1$. Then, for all $T>0$, there
exists a unique solution $f \in \mathcal{C}([0,T]; L_2^1) \cap
L^\infty(0,T; L_3^1)$ to the Boltzmann equation~(\ref{be:force})
such that $f(v,0)=f_0(v)$.
\end{theo}

\begin{proof} Let $T >0$ be fixed.
The uniqueness in $\mathcal{C}([0,T]; L_2^1) \cap L^\infty(0,T;
L_3^1)$ trivially follows from Proposition~\ref{stability}. The
proof of the existence is made in several steps, following the lines
of ~\cite[Section 3.3]{MMR}, see also \cite{MW,FoMi}.

\textit{Step 1.}  Let us first consider an initial datum $f_0 \in
L_4^1$, and  define the ``truncated'' collision operators
\begin{equation}\label{truncated}
\begin{array}{l}
\dis \IR \Q_n (f,f)(v)\, \psi(v)\d\v  =  \frac{1}{2} \IR\IR
\mathbf{1}_{\left\{|q| \leq n\right\}}\, |q|
f(v)\,f(w)\, \mathcal{A}_\zeta[\psi](v,w)\d\w\d\v, \vspace{0.2 cm}\\
\dis \IR \cL_n (f)(v)\, \psi(v)\d\v = \frac{1}{\lambda} \,\int
_{\R^3} \int _{\R^3} \mathbf{1}_{\left\{|q| \leq n\right\}}\, |q|
f(\v)\M(\w)\mathcal{J}_e[\psi](v,w) \d v \d w\,
\end{array}
\end{equation}
for any regular test function $\psi$. The operators $\Q_n$ and
$\cL_n$ are bounded in any $L_q^1$, and they are Lipschitz in
$L_2^1$ on any bounded subset of $L_2^1$. Therefore, following
\cite{Ark}, we can use the Banach fixed point theorem to get the
existence of a solution $0 \leq f_n \in \mathcal{C}([0,T]; L_2^1)
\cap L^\infty(0,T; L_4^1)$   to the Boltzmann equation $\partial_t f
=\tau \Q_n(f,f) + \cL_n(f)$. Thanks to the uniform propagation of
moments in Proposition \ref{propo:moments}, there exists a constant
$C_T>0$ (that does not depend on $n$) such that
$$
\sup_{[0,T]} \| f_n \|_{L_4^1} \leq C_T, \qquad \qquad \forall n \in \mathbb{N}.
$$

\textit{Step 2.} Let us prove that the sequence $(f_n)_n$ is
 a Cauchy sequence in $ \mathcal{C}([0,T]; L_2^1) \cap
L^\infty(0,T; L_4^1)$.  For any $m \geq n$, writing down the
equation satisfied by $f_m-f_n$ and multiplying it by
$\psi(v,t)=\mathrm{sgn}(f_m(v,t)-f_n(v,t))\langle v \rangle ^2$ as
in the proof of   Proposition \ref{stability}, we get
$$\dfrac{\d}{\d t}\int_{\R^3}|f_n(v,t)-f_m(v,t)|\langle v\rangle^2\d\v =I_{m,n}(t)+J_{m,n}(t)$$
where
$$I_{m,n}(t)=\ds \tau\int_{\R^3}
\bigg\{\Q_m(f_m,f_m)-\Q_n(f_n,f_n)\bigg\}\psi(v,t)\d\v$$  and
$$
J_{m,n}(t)=\int_{\R^3}\bigg\{\cL_m(f_m)(v,t)-\cL_n(f_n)(v,t)\bigg\}\psi(v,t)\d\v.$$
We begin by estimating $I_{m,n}(t)$. It is easy to see that
$I_{m,n}(t)=I^1_{m,n}(t)+I_{m,n}^2(t)$ where
$$I_{m,n}^1(t)=\dfrac{\tau}{2}   \IR\IR \Big( f_m(v,t)\,f_m(w,t)-f_n(v,t)f_n(
w,t)\Big)\, B_m(q) \mathcal{A}_\zeta[\psi(t)](v,w)\, \d\w\,\d\v,
$$
while $$I_{m,n}^2(t)=\ds \dfrac{\tau}{2}   \IR\IR \left(B_m(q)
-B_n(q) \right)\,f_n(v,t)f_n(w,t)
\mathcal{A}_\zeta[\psi(t)](v,w)\, \d\w\,\d\v,
$$
where $B_n(q) =|q| \mathbf{1}_{\left\{|q| \leq n\right\}}$.
Arguing as in the proof of \eqref{It}, we get easily that
$$
I_{m,n}^1(t) \leq \tau
\big \|f_n(t)+f_m(t)\big\|_{L^1_3} \big \|f_n(t)-f_m(t)\big
\|_{L^1_2}, \qquad \forall t \geq 0.
$$
The estimate of $I_{m,n}^2(t)$ is more involved. One observes
first that
$$
B_m(q)-B_n(q)=|q|\mathbf{1}_{\{n \leq |q| \leq m\}}
\leq |q| \left(\mathbf{1}_{\{|v| \geq n/2\}}+\mathbf{1}_{\{|w|\geq
n/2\}}\right).
$$
As in the proof of Proposition \ref{stability}, one has
$$
\mathcal{A}_\zeta[\psi(t)](v,w) \leq \frac{1}{4\pi}\IS \left(\langle
v'\rangle^2 + \langle w'\rangle^2 + \langle v \rangle^2 + \langle w
\rangle^2\right)\d\sigma \leq 2\left(\langle v \rangle^2 + \langle w
\rangle^2\right)
$$
and, since $|q| \leq \langle v \rangle \langle w \rangle$, one
gets
\begin{equation*}\begin{split}
I_{m,n}^2&(t) \leq \tau \IR \IR f_n(v,t)f_n(w,t)
|q|\left(\mathbf{1}_{\{|v| \geq n/2\}}+\mathbf{1}_{\{|w|\geq
n/2\}}\right)\left(\langle v \rangle^2 + \langle w
\rangle^2\right) \d\w\d\v\\
&\leq \tau \IR \IR f_n(v,t)f_n(w,t) \langle v \rangle \langle w
\rangle\left(\langle v \rangle^2 + \langle w
\rangle^2\right)\left(\mathbf{1}_{\{\langle v \rangle \geq
n/2\}}+\mathbf{1}_{\{\langle w \rangle\geq n/2\}}\right) \d\w\d\v.
\end{split}\end{equation*}
It is not difficult to deduce then that
$$
I_{m,n}^2(t) \leq 4\tau \bigg(\IR f_n(v,t) \langle \v \rangle^3 \d\v \bigg) \bigg(\IR
f_n(v,t) \langle \v \rangle^3 \mathbf{1}_{\{\langle \v \rangle \geq
n/2\}}\d\v\bigg).
$$
Since $\sup_{[0,T]} \|f_n(t)\|_{L^1_4} \leq C_T$ for any $n \in
\mathbb{N},$ the latter integral is estimated as
$$
\IR
f_n(v,t) \langle \v \rangle^3 \mathbf{1}_{\{\langle \v \rangle \geq
n/2\}}\d\v \leq \IR f_n(v,t) \langle \v \rangle^4
\mathbf{1}_{\{\langle \v \rangle \geq n/2\}}\frac{\d\v}{\langle v
\rangle} \leq \frac{2C_T}{n}
$$
and we get
$$
I_{m,n}^2(t)\leq 4\tau \bigg(\IR f_m(v,t) \langle \v \rangle^3 \d\v
\bigg)\dfrac{2C_T}{n} \leq \dfrac{8C^2_T\tau}{n}, \qquad \forall t
\in [0,T], \quad m \geq n.
$$
Therefore,
\begin{equation}
I_{m,n}(t) \leq 2\tau C_T \|f_n(t)-f_m(t)\big \|_{L^1_2} +
\dfrac{8C^2_T\tau}{n}, \qquad \forall t \in [0,T], \quad m \geq n.
\end{equation}
We proceed in the same way with $J_{m,n}(t)$. First, we notice that
$J_{m,n}(t)$ splits as $J_{m,n}(t)=J_{m,n}^1(t)+J_{m,n}^2(t)$ with
$$
J_{m,n}^1(t)=\dfrac{1}{\lambda}\IR \IR
B_m(q)\left[f_m(v,t)-f_n(v,t)\right
]\M(w)\mathcal{J}_e\left[\psi(t)\right](v,w)\d\v\d\w
$$
and
$$
J_{m,n}^2(t)=\dfrac{1}{\lambda}\IR \IR \left[B_m(q)-B_n(q)\right]
f_n(v,t) \M(w)\mathcal{J}_e\left[\psi(t)\right](v,w)\d\v\d\w.
$$
Arguing as in the proof of Proposition \ref{stability}, we get
$$J_{m,n}^1(t) \leq \IR\!\! \cL_{m}(|f_n-f_m|)(v,t) \langle \v \rangle^2 \d\v
\leq \frac{2\kappa}{\lambda}\IR\! \IR\! |q|^2|w||(f_n-f_m)(v,t)|\M(w) \d\v\d\w
$$
and there exists a
positive constant $c_+$ such that
$$
J_{m,n}^1(t) \leq c_+ \big\|f_n(t)-f_m(t)\|_{L^1_2}, \qquad
\forall t \in [0,T].
$$
Let us now estimate $J_{m,n}^2(t)$. As above,
$$
J_{m,n}^2(t)= \dfrac{1}{\lambda}\IR \IR |q| \mathbf{1}_{\{n \leq |q|
\leq m\}}f_n(v,t)
\M(w)\mathcal{J}_e\left[\psi(t)\right](v,w)\d\v\d\w
$$
and
$$
\mathcal{J}_e\left[\psi(t)\right](v,w) \leq
\dfrac{1}{2\pi}\int_{\S} |\hat{q}\cdot\n|\left(\langle
v^\star\rangle^2+\langle v \rangle^2\right)\d\n
=\mathcal{J}_e\left[\langle \cdot \rangle^2\right] +
\frac{2}{2\pi}\int_{\S} |\hat{q}\cdot\n| \langle v \rangle^2 \d\n.
$$
Calculations already performed lead then to
$$
\mathcal{J}_e\left[\psi(t)\right] (v,w) \leq -2\kappa \langle
q,w\rangle + 2\langle v\rangle ^2 \leq 2\left(\langle v \rangle
\langle w \rangle^2 + \langle v \rangle ^2 \right), \qquad \forall\,
v,w \in \R^3\,.
$$
Finally,
\begin{equation*}\begin{split}
J_{m,n}^2(t) &\leq \frac{2}{\lambda} \IR \IR |q|\mathbf{1}_{\{|q|
\geq n\}}f_n(v,t)\M(w)\left(\langle v \rangle \langle w \rangle^2 +
\langle v \rangle ^2 \right)\d\w\d\v\\
&\leq \frac{2}{\lambda} \IR \IR \mathbf{1}_{\{|q| \geq
n\}}f_n(v,t)\M (w)\left(\langle v \rangle^3 \langle w \rangle  +
\langle v \rangle ^2 \langle w \rangle^3 \right)\d\w\d\v.
\end{split}\end{equation*}
Now, arguing as we did for $I_{m,n}^2(t)$, there exists some
constant $\tilde{C}_T$ that depends only on $\|\M\|_{L^1_4}$ and
$\sup_{n}\sup_{[0,T]}\|f_n(t)\|_{L^1_4}$  such that
$$
J_{m,n}^2(t) \leq \dfrac{\tilde{C}_T}{n}, \qquad \forall m \geq n,\quad t \in [0,T].
$$
Gathering all these estimates, we obtain the existence of constants
$C_1(T)$ and $C_2(T)$ that do not depend on $m,n$ such that
$$
\dfrac{\d}{\d t}\int_{\R^3}|f_n(v,t)-f_m(v,t)|\langle
v\rangle^2\d\v \leq C_1(T) \|f_n(t)-f_m(t)\big \|_{L^1_2} +
\dfrac{C_2(T)}{n}, \qquad \forall t \in [0,T], \: m \geq n.
$$

This ensures that $(f_n)_n$ is a Cauchy sequence in
$\mathcal{C}([0,T]; L_2^1)$. Denoting by $f $ its limit, we obtain
that $f\in \mathcal{C}([0,T]; L_2^1) \cap L^\infty(0,T; L_4^1)$ is a
solution to the Boltzmann equation~(\ref{be:force}) (with the actual
collision operators $\Q$
and $\cL$).\\

\textit{Step 3.} When the initial datum $f_0 \in L_3^1$, we
introduce the sequence of initial data $f_{0,j} := f_0
\mathbf{1}_{|v| \leq j}$. Since $f_{0,j} \in L_4^1$, we have the
existence of a solution $f_j \in \mathcal{C}([0,T]; L_2^1) \cap
L^\infty(0,T; L_4^1)$ to the Boltzmann equation associated to the
initial datum $f_{0,j}$. Moreover, there exists $C_T$ such that
$\sup_{[0,T]} \| f_j \|_{L_3^1} \leq C_T$. We establish again that
$(f_j)_j$ is a Cauchy sequence in $\mathcal{C}([0,T]; L_2^1)$ by
using the $L^1$-stability in Proposition \ref{stability}.
\end{proof}

\section{Existence of non-trivial stationary state}

All the material of the previous sections allows us to state our
main result:
\begin{theo}[\textit{\textbf{Existence of stationary solutions}}] \label{theo:main}
For any distribution function $\M(\v)$ satisfying Assumption
\ref{AssF1} and  any $\tau \geq 0$, there exists a nonnegative $F
\in L^1_2\cap L^p$, $p\in (1,\infty)$ with unit mass and positive
temperature such that $\tau \Q(F,F) +\cL(F)=0$.
\end{theo}

\begin{proof} As already announced, the existence of stationary solution
to \eqref{be:force} relies on the application of Lemma \ref{GPV}
to the evolution semi-group $(\mathcal{S}_t)_{t \geq 0}$ governing
\eqref{be:force}. Namely, for $f_0 \in L^1$, let
$f(t)=\mathcal{S}_t f_0$ denote the unique solution to
\eqref{be:force} with initial state $f(0)=f_0$.  The continuity
properties of the semi-group are proved by the study of the Cauchy
problem, recalled in Section~\ref{sec:existence}. Let us fix
$p_o\in(1,\infty)$. On the Banach space $\mathcal{Y}= L^1 _2$,
thanks to the uniform bounds on the $L^1 _3$ and $L^{p_o}$ norms,
the nonempty convex subset
     \[ \mathcal{Z}= \left\{ 0 \le f \in \mathcal{Y}, \hspace{0.3cm}
         \int_{\R^3} f  \, \d v = 1  \quad \mbox{ and }
         \quad \|f\|_{L^1 _3} + \|f\|_{L^{p_o}} \le M \right\} \]
is stable by the semi-group provided $M$ is big enough. This set
is weakly compact in $\mathcal{Y}$ by Dunford-Pettis Theorem, and
the continuity of $\mathcal{S}_t$ for all $t \ge 0$ on
$\mathcal{Z}$ follows from Proposition~\ref{stability}. Then,
Lemma \ref{GPV} shows that there exists a nonnegative stationary
solution to \eqref{be:force} in $L^1 _3 \cap L^{p_o}$ with unit
mass. In fact, the uniform in time $L^p$ bounds also imply the
boundedness of $F$ in $L^p$ for all $p\in (1,\infty)$.
\end{proof}

 As a
corollary of Theorem \ref{theo:main}, choosing $\tau=0$ allows us to
prove the existence of a steady state to the \textbf{\textit{linear
inelastic scattering operator}} $\cL$ when the distribution function
of the background is not a Maxwellian, generalizing the result of
\cite{MaPi,LoTo,SpTo}.

\begin{cor} Let $\M$ satisfy Assumption \ref{AssF1}. Then, the linear inelastic scattering operator
$\cL$ defined by \eqref{linearoperator} admits a unique nonnegative
steady state $F \in L^1_2\cap L^p$, $p\in (1,\infty)$, with unit
mass and positive temperature.
\end{cor}

\begin{proof} The existence of a nonnegative equilibrium solution $F \in L^1_2$
is a direct application of Theorem \ref{theo:main} with $\tau=0.$ A
simple application of  Krein-Rutman  Theorem implies the uniqueness
of the stationary solution $F$ within the range of nonnegative
distributions with unit mass.
\end{proof}
\begin{nb}[\textit{\textbf{$H$-Theorem and trend towards equilibrium}}]\label{Htheo}
As in \cite{LoTo}, it is possible to prove a linear version of the
classical $H$-Theorem for the linear inelastic Boltzmann equation
\eqref{be:force} with $\tau=0:$
\begin{equation}
\label{be:forcelinear} \partial_t f  = \mathcal{L}(f), \qquad
f(t=0)=f_0 \in L^1.
\end{equation} Namely,  for any {\it convex}
$\mathcal{C}^1$ function $\Phi \::\:\mathbb{R}^+ \to \mathbb{R}$,
let
$$H_{\Phi}(f|F)=\int_{\R^3}F(\v) \Phi\left (\dfrac{f(\v)}{F(\v)}
\right)\d\v, \qquad f \in L^1.$$ Arguing as in \cite{LoTo}, it is
easy to prove that, if the initial state $f_0$ has unique mass and
finite entropy $H_{\Phi}(f_0|F) < \infty$, then
\begin{equation}\label{h}
\dfrac{\d}{\d t}H_{\Phi}(f(t)|F)\leq 0 \qquad \qquad (t \geq 0)
\end{equation}
where $f(t)$ stands for the (unique) solution to
\eqref{be:forcelinear}. Moreover, still arguing as in \cite{LoTo},
one proves that if moreover $\dis \IR \left(1 + \v^2 + |\log
f_0(\v)|\right) f_0(\v) \d\v < \infty$, then
$$\lim_{t \to \infty}\IR \left|f(v,t)-F(v)\right|\d\v=0.$$
\end{nb}

\section{Regularity of the steady state}

In this final section, our aim is to establish the existence of some
smooth stationary solution to~\eqref{be:force}. Namely, adopting the
strategy of \cite[Section 4.1]{MiMo}, we  prove
\begin{theo}[\textit{\textbf{Regularity of stationary solutions}}]  \label{regul} There exists a stationary solution $F$ to the
Boltzmann equation
$$\tau \Q(F,F) +\cL(F)=0
$$
that belongs to $C^\infty(\R^3)$.
\end{theo}
We shall follow the same lines of \cite[Theorem 5.5]{MoVi} and \cite[Section 3.6]{MiMo}, from which we deduce the
exponential decay in time of singularities and thus the smoothness
of stationary solutions. This proof needs the following ingredients:
\begin{enumerate}[i)]
\item The stability result already proved in
Proposition~\ref{stability}.
\item An estimate on the Duhamel
representation \cite[Proposition 3.4]{MiMo} of the solution to
\eqref{be:force} (see Proposition \ref{Duhamel}).
\item A result of propagation of Sobolev norms (see Proposition
\ref{SobNorms}).
\end{enumerate}

Let us first extend the regularity estimate of \cite[Proposition
3.4]{MiMo} to our situation. For any $f \in L^1$, let
$$\Sigma(f)(v)=\tau \left(|\cdot| * f\right)(v) + \nu(v) =\tau \IR |v-w|f(w)\d\w +
\nu(v).$$ It is easy to see that, for $f_0 \in L^1_3$, the unique
solution $f(t)$ to \eqref{be:force} is given by the the following
Duhamel representation:
\begin{equation}\begin{split} \label{Duhrepr} f(v,t) &=\ds f_0(v)\, {\rm e}^{-\int_0^t \Sigma(f)(v,s) \d s} + \int_0^t \big( \tau\,\Q^+ (f,f) +
\cL^+(f) \big) (v,s)\, {\rm e}^{-\int_s^t \Sigma(f)(v,r) \d r }\d s\\
&=f_0(v)G(v, 0,t)+ \int_0^t  \left(\tau\, \Q^+(f,f)+\cL^+(f)
\right)(v,s)G(v,s,t)\d s
\end{split}
\end{equation}
where we set
$$G(v,s,t) =\exp\left(-\int_s^t \Sigma(f)(v,r) \d r \right)\,\qquad 0 \leq
s \leq t, \qquad v \in \R^3.$$
\begin{propo} \label{Duhamel} There are some positive constants  $C_\mathrm{Duh}, K$ such
that for any $k\in \mathbb{N}$ and $\eta \geq 0$ we have
\begin{equation}\label{eq.5.2}
\| f_0(\cdot)G(\cdot,0,t) \|_{H_\eta^{k+1}} \leq C_{\rm Duh}\,
{\rm e}^{-Kt} \| f_0 \|_{H_{\eta+1}^{k+1}} \!\left(\sup_{0 \leq r
\leq t}\| f(\cdot, r)\|_{H_{\eta+3}^{k}}^2 \! \! \!+  \! \!\sup_{0
\leq r \leq t}\| f(\cdot, r)\|_{H_{\eta+3}^{k}}^{k+3}\right)
\end{equation}
and
\begin{align}\label{eq.5.3}
&\left\| \int_0^t G(\cdot,s,t)\big( \tau\, \Q^+(f,f) + \cL^+(f)
\big) (\cdot,s) \d s \right\|_{H_\eta^{k+1}} \\
&\qquad \qquad\qquad\qquad\qquad\qquad\leq C_{\rm Duh}\,
\left(\sup_{0 \leq r \leq t}\| f(\cdot, r)\|_{H_{\eta+3}^{k}}^2 +
\sup_{0 \leq r \leq t}\| f(\cdot,
r)\|_{H_{\eta+3}^{k}}^{k+3}\right) .\nonumber
\end{align}
\end{propo}

\begin{proof}
The proof is quite similar to \cite[Proposition 5.2]{MoVi}. Here,
for simplicity we have done it for natural $k$, although it is
simple to generalize it to $k>0$ by interpolation. Precisely, for
any $f \in L^1 $ define
$$
L(f)(\v)=\int_{\R^3} |\v-\w|f(\w)\d\w.
$$
It is clear that
$$
\Sigma(f)(\v)=L(\tau f+\M)(\v), \qquad \forall f \in L^1.
$$
Now, according to \cite[Lemma 4.3]{duduch}, for any given $k \geq 0$
and any $\delta > 3/2$, the linear operator
$$
L \::\:H^k_\delta \longrightarrow W_{-1}^{k+1,\infty}
$$
is bounded, i.e. for any $\delta >3/2$ and any $k \geq 0$, there
exists $C_{k,\delta}$ such that
$$
\left\|L(g)\right\|_{W_{-1}^{k+1,\infty}} \leq C_{k,\delta}
\|g\|_{H^k_\delta}, \qquad \forall g \in H_\delta^k.
$$
Let us fix now $ k \in \mathbb{N}$ and $\delta >3/2$. Since $\M
\in H_\delta^k$ due to Assumption \ref{AssF1}, one deduces that
$$
\left\|\Sigma(f)\right\|_{W_{-1}^{k+1,\infty}} \leq C
\|f+ \M\|_{H^k_\delta}, \qquad \forall f \in H_\delta^k
$$
where, as in the rest of the proof, we shall denote any positive
constant independent of $f$ and possibly dependent on $\M$ by $C$.
Setting
$$
F(v,s,t)=\ds \int_s^t \Sigma(f)(v,r)\d r,
$$
one sees that
$$
\left\|F(\cdot,s,t)\right\|_{W_{-1}^{k+1,\infty}} \leq C
\sqrt{t-s}\,\bigg(\int_s^t \|f(\cdot,r)\|_{H^k_\delta}^2\d r
\bigg)^{1/2}+C (t-s)\, \|\M\|_{H^k_\delta}, \qquad 0 \leq s \leq
t.
$$
Now, since $L(g) \geq 0$ for any $g \geq 0$, according to
Assumption \ref{AssF1} and \eqref{hypp}, we see that there exists
some constant $\chi>0$ such that
$$
\Sigma(f)(v) \geq L(\M)(\v) \geq \chi, \qquad \forall f \in L^1,\:f
\geq 0,\;\forall v \in \R^3.
$$
By taking the successive derivatives of
$G(v,s,t)=\exp(-F(v,s,t))$, one gets as in \cite[Proposition
5.2]{MoVi}
\begin{equation}\begin{split} \label{estG}
\| G( \cdot,s,t) \|_{W_{-1}^{k+1,\infty}} &\leq  C{\rm
e}^{-\chi(t-s)} \left[\sqrt{t-s}\,\bigg(\int_s^t
\|f(\cdot,r)\|_{H^k_\delta}^2\d r
\bigg)^{(k+1)/2}\!\!\!\!\!\!+(t-s)\, \|\M\|_{H^k_\delta}+1\right]\\
&\leq C {\rm e}^{-K(t-s)} \left( 1+ \sup_{s \leq r \leq t} \|
f(\cdot, r) \|_{H_\delta^k}^{k+1}\right) ,\end{split}
\end{equation}
for some $0 < K < \chi$. Then, we shall use the following estimate
\cite[Lemma 5.3]{MoVi} that allows to exchange a time integral and
a Sobolev norm:
$$
\left\| \int_0^t Z(\cdot,s) \d s \right\|_{H_\ell^r} \leq
\frac{1}{\sqrt{\lambda}} \left( \int_0^t {\rm e}^{\lambda (t-s)} \|
Z(\cdot,s) \|_{H_\ell^r}^2 \d s \right)^{1/2}, \qquad \forall
\lambda>0,\,\,\, \forall \ell, r \in \R\,.
$$
As a consequence we have for any $k \geq 0$,
$$
\begin{array}{c}
\dis \left\| \int_0^t \big( \tau\,\Q^+(f,f) + \cL^+(f) \big) (\cdot,s)\,\, G(\cdot,s,t) \d s \right\|_{H_\eta^{k+1}} \vspace{0.2 cm}\\
\dis \leq C \left[ \int_0^t {\rm e}^{K(t-s)} \left( \left\|
\tau\,\Q^+(f,f)(\cdot,s) \right\|_{H_{\eta+1}^{k+1}}^2 + \left\|
\cL^+(f)(\cdot,s) \right\|_{H_{\eta+1}^{k+1}}^2 \right) \left\| G(
\cdot,s,t) \right\|_{W_{-1}^{k+1,\infty}}^2\,\, \d s \right]^{1/2}.
\end{array}
$$
Recall now the  so--called Bouchut-Desvillettes-Lu regularity result
in Propositions \ref{Q^+} and \ref{L^+}:
$$
\| \Q^+(f,f) \|_{H_{\eta+1}^{k+1}} \leq C \Big[  \| f
\|_{H_{\eta+3}^{k}}^2 +  \| f \|_{L_{\eta+3}^1}^2 \Big]
$$
and
$$
     \left\|\cL^+(f)\right\|_{H^{k+1}_{\eta+1}} \le
     C \left[\left\|\M\right\|_{H^{k}_{\eta+3}} \left\|f\right\|_{H^{k}_{\eta+3}} +
     \left\|\M\right\|_{L^1_{\eta+3}}
     \left\|f\right\|_{L^1_{\eta+3}}\right].
$$
Arguing now as in \cite[Proposition 5.2]{MoVi} and using the
estimate~(\ref{estG}) with the choice $\delta=\eta+3$, we get
\begin{equation*}\begin{split}
&\left\| \int_0^t \big( \tau\,\Q^+(f,f) + \cL^+(f) \big) (\cdot,s)\,\, G(\cdot,s,t) \d s \right\|_{H_\eta^{k+1}} \vspace{0.2 cm}\\
&\leq C \left[ \int_0^t {\rm e}^{K(t-s)} \left\| f(\cdot,s)
\right\|_{H_{\eta+3}^{k }}^4  {\rm e}^{-2K(t-s)} \left(1+\sup_{s
\leq r \leq t} \| f(\cdot, r)\|_{H_{\eta+3}^k}^{k+1}\right)^2\,\,
\d s
\right]^{1/2}\\
&\dis \leq C \max \left(\sup_{0 \leq r \leq t}\| f(\cdot,
r)\|_{H_{\eta+3}^{k}}^2,\sup_{0 \leq r \leq t}\| f(\cdot,
r)\|_{H_{\eta+3}^{k}}^{k+3}\right)
\end{split}\end{equation*}
which proves \eqref{eq.5.3}. The proof of \eqref{eq.5.2} is similar.
\end{proof}

A direct consequence of the previous result together with the
uniform $L^2$ bounds is the uniform in time propagation of Sobolev
norms. The proof is carried on exactly as in
\cite[Proposition~3.5]{MiMo}.

\begin{propo} \label{SobNorms} Let $\M$ satisfy Assumption
\ref{AssF1}. Let $f_0\in L^1_2$, $f_0 \geq 0$ with unit mass and let
$f$ be the unique solution of the Boltzmann
equation~\eqref{be:force} in $\mathcal{C}(\mathbb{R}^+; L_2^1) \cap
L^\infty(\mathbb{R}^+; L_3^1)$ associated with~$f_0$. Then, for all
$s>0$ and $\eta \geq 1$, there exists $w(s)>0$ such that
$$
f_0 \in H_{\eta+w}^s \quad \Longrightarrow \quad \sup_{t \geq 0}
\| f(\cdot,t) \|_{H_\eta^s} < +\infty\,.
$$
\end{propo}

The previous ingredients allow to proof the following theorem, see
\cite[Theorem 5.5]{MoVi} for the proof.

\begin{theo}[\textit{\textbf{Exponential decay of singularities}}] \label{decaysing}
Let $f_0\in L^1_2\cap L^2$ with unit mass and let $f$ be the
unique solution of the Boltzmann equation~\eqref{be:force} in
$\mathcal{C}(\mathbb{R}^+; L_2^1) \cap L^\infty(\mathbb{R}^+;
L_3^1)$ associated with~$f_0$. Let $\M$ satisfy Assumption
\ref{AssF1}.  Let $s \geq 0$, $q \geq 0$ be arbitrarily large.
Then $f$ splits into the sum of a regular and a
singular part $f=f_{R} +f_{S}$ where
$$
\left\{
\begin{array}{l}
\dis \sup_{t \geq 0} \| f_R(t) \|_{H_q^s \cap L_2^1} < + \infty\,, \qquad f_R \geq 0 \vspace{0.2 cm}\\
\dis \exists \lambda>0\,: \qquad \| f_S(t) \|_{L_2^1} = O \big( {\rm
e}^{-\lambda t} \big)\,.
\end{array}
\right.
$$
\end{theo}

\begin{proof} The proof is easily adapted from \cite[Theorem 3.6]{MiMo} since the
$L^1$-stability result (Proposition~\ref{stability}), the Duhamel
representation (Proposition~\ref{Duhamel}), the uniform propagation of
Sobolev norms (Proposition~\ref{SobNorms}) allow to adapt directly
\cite[Theorem 5.5]{MoVi}.
\end{proof}

Finally, Theorem \ref{decaysing} allows to prove the main
Theorem~\ref{regul}.
 \bigskip

\noindent {\bf Acknowledgements:} JAC acknowledges the support
from DGI-MEC (Spain) FEDER-project MTM2005-08024. JAC and MB
acknowledge partial support of the Acc. Integ. program
HI2006-0111. MB acknowledges support also from MIUR (Project
``Non--conservative binary interactions in various types of
kinetic models''), from GNFM--INdAM, and from the University of
Parma. JAC acknowledges IPAM--UCLA where part of this work was
done. Finally, BL wishes to thank Cl\'{e}ment Mouhot for fruitful
discussions.

\end{document}